\documentclass[11pt]{article}

\usepackage{amssymb}
\usepackage{amsmath}
\usepackage{amsfonts}
\usepackage{epsfig}

\newcommand{\R}{\mathbb{R}}
\newcommand{\I}{\mathbb{I}}
\newcommand{\E}{\mathbb{E}}

\newcommand{\Id}{\mbox{Id}}
\newcommand{\CQFD}
{%
\mbox{}%
\nolinebreak%
\hfill%
\rule{2mm}{2mm}%
\medbreak%
\par%
}

\newtheorem{proposition}{Proposition}
\newtheorem{corollary}{Corollary}
\newtheorem{lemma}{Lemma}

\begin{document}

\title{Gaussian Regularized Sliced Inverse Regression}
\author{Caroline Bernard-Michel, Laurent Gardes, St\'ephane Girard$^{\star}$}
\date{Laboratoire Jean-Kuntzmann \& INRIA Rh\^one-Alpes, team Mistis,\\
 Inovall\'ee, 655, av. de l'Europe, Montbonnot, 38334 Saint-Ismier cedex, France, \\
($^{\star}$ corresponding author, {\tt Stephane.Girard@inrialpes.fr})
}
\maketitle
\begin{abstract}
Sliced Inverse Regression (SIR) is an effective method for dimension
reduction in high-dimensional regression problems. The original
method, however, requires the inversion of the predictors covariance
matrix. In case of collinearity between these predictors or small sample
sizes compared to the dimension, the inversion is not possible
and a regularization technique has to be used.
Our approach is based on a Fisher Lecture given by R.D. Cook
where it is shown that SIR axes can be
interpreted as solutions of an inverse regression problem. 
In this paper, a Gaussian prior distribution is introduced on the unknown
parameters of the inverse regression problem in order
to regularize their estimation. We show that some existing SIR
regularizations can enter our framework,
which permits a global understanding of these methods.
Three new priors are proposed leading to new regularizations
of the SIR method. A comparison on simulated data is provided.

\noindent {\bf Keywords:} Inverse regression, regularization, sufficient dimension reduction.  

\end{abstract}
%%%%%%%%%%%%%%%%%%%%%%%%%%%%%%%%%%%%%%%%%%%%%%%%%%%%%%%%%%%%%%%%%%%%%%%
\section{Introduction}
%%%%%%%%%%%%%%%%%%%%%%%%%%%%%%%%%%%%%%%%%%%%%%%%%%%%%%%%%%%%%%%%%%%%%%%

Many methods have been developed for inferring the conditional distribution
of an univariate response $Y$ given a predictor $X$ in $\R^p$,
ranging from linear regression~\cite{Regression}
to support vector regression~\cite{CST}.
When $p$ is large, sufficient dimension reduction aims at replacing
the predictor $X$ by its projection onto a subspace of smaller dimension
without loss of information on the conditional distribution of
$Y$ given $X$. In this context, the central subspace, denoted by
$S_{Y|X}$ plays an important role. It is defined as the smallest
subspace such that, conditionally on the projection
of $X$ on $S_{Y|X}$, $Y$ and $X$ are independent. In other words, the projection
of $X$ on $S_{Y|X}$ contains all the information on $Y$ that is available in 
the predictor $X$.

The estimation of the central subspace has received considerable
attention since the past twenty years.  Without intending to be
exhaustive, we refer to
Sliced Inverse Regression (SIR)~\cite{Li}, sliced average
variance estimation~\cite{CookWeis},
and graphical regression~\cite{Cook98} methods.
Among them, SIR seems to be the most popular one. The original
method has been adapted to various frameworks and the relative
asymptotic properties have been derived, see for
instance~\cite{HsingCarr, ZhuNg, Saracco97, GaSa03, Saracco05}.
We also refer to~\cite{AraSara} for a study of the SIR finite
sample properties, to~\cite{Ferre,Schott} for the estimation of 
$K=\dim(S_{Y|X})$, the dimension of the central subspace,
and to \cite{Anes,FerYao} for extension to functional covariates.

Assuming that $K$ is known and introducing $\Sigma=\mbox{Cov}(X)$,
in the SIR methodology,
a basis of the central subspace is obtained by computing the eigenvectors
associated to the largest $K$ eigenvalues of $\Sigma^{-1}\mbox{Cov}(\E(X|Y))$.
Unfortunately, the classical $n-$ sample estimate $\hat\Sigma$ of $\Sigma$ can
be singular, or at least ill-conditioned, in several situations.
Indeed, since rank$(\hat\Sigma)\leq\min(n-1,p)$, if $n\leq p$ then $\hat\Sigma$
is singular. Even when $n$ and $p$ are of the same order,
$\hat\Sigma$ is ill-conditioned, and its inversion introduces
numerical instabilities in the estimation of the central subspace.
Similar phenomena occur when the coordinates of $X$ are highly
correlated. 

Some regularizations of the SIR method have been proposed
to overcome this limitation.
In~\cite{Chiar} and~\cite{LiLi}, a Principal Component Analysis (PCA)
is used as a preprocessing step in order to eliminate the directions
in which the random vector $X$ is degenerated. Thus, for a properly
chosen dimension $d$ of the projection subspace, the covariance
matrix of the projected observations is regular.
In the sequel, this technique will be referred to as PCA+SIR.
Another method consists in adopting a ridge regression technique
(see for instance~\cite{Regression}, Chapter~17) {\it i.e.}
replaces the sample estimate $\hat\Sigma$
by a perturbed version $\hat\Sigma+\tau I_p$ where $I_p$ is the $p\times p$
identity matrix and $\tau$ is a positive real number~\cite{Zhong}.
Here, the idea is that, for $\tau$ large enough, $\hat\Sigma+\tau I_p$
is regular and its condition number increases with $\tau$.
Similarly, in~\cite{Scrucca,Scrucca2},
regularized discriminant analysis~\cite{Fried}
is adapted to the SIR framework.
More recently, it is proposed in~\cite{LiYin} to interpret SIR 
as an optimization problem and to introduce $L_1-$ and $L_2-$
penalty terms in the optimized criterion.

Our approach is based on a Fisher Lecture given by R.D. Cook~\cite{Cook07}
where it is shown that the axes spanning the central subspace can be
interpreted as the solutions of an inverse regression problem. 
In this paper, a Gaussian prior is introduced on the unknown
parameters of the inverse regression problem in order
to regularize their estimation. We show that the previously
mentioned techniques~\cite{Chiar,LiLi,Zhong} can enter our framework,
which permits a global understanding of these methods.
Three new priors are proposed leading to new regularizations
of the SIR method. A comparison with the $L_2-$ regularization
introduced in~\cite{LiYin} is
also provided. It is shown that, from the theoretical point of view,
 the proposed $L_2-$ regularization cannot be justified.

This paper is organized as follows. In Section~\ref{IRWR},
an adaptation of the inverse regression model to our framework
is presented. Section~\ref{RIR} is dedicated to the regularization
aspects. Theoretical comparisons with existing approaches as well as new 
methods are provided. Finite sample properties are illustrated
in Section~\ref{FSP}. Proofs are postponed to the Appendix.

%%%%%%%%%%%%%%%%%%%%%%%%%%%%%%%%%%%%%%%%%%%%%%%%%%%%%%%%%%%%%%%%%%%%%%%
\section{Inverse regression without regularization}
%%%%%%%%%%%%%%%%%%%%%%%%%%%%%%%%%%%%%%%%%%%%%%%%%%%%%%%%%%%%%%%%%%%%%%%
\label{IRWR}

Consider $X$ a $\R^p-$ random vector, $Y$ the real response variable
and let us denote by $S_{Y|X}$ the central subspace.
In the following, for the sake of simplicity, we assume that
$K=\dim(S_{Y|X})=1$, the case $1<K<p$ being discussed in 
Section~\ref{RIR}. We thus introduce $\beta\in\R^p$ such
that $S_{Y|X}=\mbox{span}(\beta)$.
In Subsection~\ref{parIRM}, 
the considered inverse regression model is presented.
The estimation of the unknown parameters is discussed in
Subsection~\ref{MLE} and the links with the SIR method
are established in Subsection~\ref{SIR}.

\subsection{Single-index inverse regression model}
%%%%%%%%%%%%%%%%%%%%%%%%%%%%%%%%%%%%%%%%%%%%%%%%%%%%%%%%%%%%%%%%%%%%%%%
\label{parIRM}

The following inverse single-index regression
model is considered (see~\cite{Cook07}, equation~(2) for the
multi-index model):
\begin{equation}
\label{IRM}
X=\mu + c(Y) V  b +\varepsilon,
\end{equation}
where $\mu$ and $b$ are non-random $\R^p-$ vectors, $\varepsilon$ 
is a centered $\R^p-$ Gaussian random vector, independent of $Y$,
with covariance matrix $\mbox{Cov}(\varepsilon)=V$
and $c:\R\to\R$ is a nonrandom function. 
Under this model, 
$$
\E(X|Y=y)=\mu + c(y) V  b,
$$
and thus, after translation by $\mu$, the conditional expectation of $X$
given $Y$ is a degenerated random vector located in the direction $V b$. 
From~\cite{Cook07}, Proposition~1, $b$ corresponds to the direction $\beta$
of the central subspace.
In the sequel, it will appear that, under appropriate conditions,
the maximum likelihood estimator of $b$ is (up to a scale parameter)
the SIR estimator of $\beta$.
Moreover, note that, under~(\ref{IRM}), one has 
\begin{equation}
\label{eqc}
c(y)=\frac{\E(b^t (X-\mu) |Y=y)}{b^t V b}.
\end{equation}
Now, restricting ourselves to the single-index case, the forward model of SIR
asserts that there exists a univariate link function $g$ such that 
$\E(Y|X)=g(b^t X)$ or equivalently, $b^tX=g^{-1}(\E(Y|X))$.
Thus, replacing in~(\ref{eqc}) yields 
$$
c(y)=\frac{g^{-1}(y)-b^t\mu}{b^t V b},
$$
{\it i.e.} the coordinate function is, up to an affine function,
the inverse of the link function in the single-index forward model of SIR.

\subsection{Maximum likelihood estimation}
%%%%%%%%%%%%%%%%%%%%%%%%%%%%%%%%%%%%%%%%%%%%%%%%%%%%%%%%%%%%%%%%%%%%%%%
\label{MLE}

We now address the estimation of the coordinate function $c(.)$,
the direction $b$, the covariance matrix $V$ and the location parameter $\mu$
in model~(\ref{IRM}).
To this end, we focus
on a projection estimator of the unknown function $c(.)$. More
precisely, it is expanded as a linear combination of $h$ basis
functions $s_j(.)$, $j=1,\dots,h$:
\begin{equation}
\label{projc}
c(.)=\sum_{j=1}^h c_j s_j(.),
\end{equation}
where the coefficients $c_j$, $j=1,\dots,h$ are unknown whereas
$h$ is supposed to be known.
Introducing $c=(c_1,\dots,c_h)^t$ and $s(.)=(s_1(.),\dots,s_h(.))^t$,
model~(\ref{IRM}) can be rewritten as
\begin{equation}
\label{IRMproj}
X=\mu +   s^t(Y) c V b +\varepsilon,\;\; \varepsilon\sim  {\cal N}(0,V),
\end{equation}
where ${\cal N}(0,V)$ is the multivariate centered Gaussian
distribution with covariance matrix~$V$.
In the sequel we denote by
$$
\rho = \frac{ b^t \Sigma b - b^t V b}{b^t V b},
$$
the signal to noise ratio in the direction $b$.
Let $(X_i,Y_i)$, $i=1,\dots,n$ be a sample of independent random variables
distributed as $(X,Y)$.
Clearly, estimating $(\mu,V,b,c)$ by maximization of the likelihood 
in model (\ref{IRMproj}) consists in minimizing
\begin{equation}
\label{Gproj}
G(\mu,V,b,c)= \log\det V +
 \frac{1}{n}\sum_{i=1}^n (\mu+ s^t(Y_i) c Vb-{X}_i )^t V^{-1} 
(\mu+ s^t(Y_i) c Vb-{X}_i ),
\end{equation}
with respect to $(\mu,V,b,c)$. 
Note that $G(\mu,V,b,c)$ can also
be interpreted as a discrepancy functional, see equation~(5) in~\cite{CookNi}.
Up to our knowledge, the introduction of
such functional in the inverse regression framework is due to~\cite{Cook04}.
Let us introduce the
$h\times h$ empirical covariance matrix $W$ of $s(Y)$ defined by
$$
W=\frac{1}{n}\sum_{i=1}^n (s(Y_i)-\bar{s}) (s(Y_i)-\bar{s})^t,
$$
the $h\times p$ matrix $M$ defined by
$$
M=\frac{1}{n}\sum_{i=1}^n (s(Y_i)-\bar{s}) (X_i-\bar{X})^t,
$$
and $\hat\Sigma$ the empirical $p\times p$ covariance matrix of $X$
$$
\hat\Sigma=\frac{1}{n}\sum_{i=1}^n (X_i-\bar{X}) (X_i-\bar{X})^t,
$$
where 
$$
\bar{X}=\frac{1}{n}\sum_{i=1}^n X_i\; \mbox{ and }\;
\bar{s}=\frac{1}{n}\sum_{i=1}^n s(Y_i).
$$
\begin{lemma}
\label{lemmaG}
Using the above notations, 
$G(\mu,V,b,c)$ can be rewritten as
\begin{eqnarray*}
{G}(\mu,V,b,c)&=& \log\det V + \mbox{tr}(\hat \Sigma V^{-1}) 
+ (\mu-\bar{X}+\bar{s}^tc Vb)^t V^{-1}(\mu-\bar{X}+\bar{s}^tc Vb)\\
&+& (c^t W c) (b^t V b) -  2 c^t M b.
\end{eqnarray*}
\end{lemma}
The maximum likelihood estimators of $\mu$, $V$, $b$ and
$c$ are closed-form.
\begin{proposition}
\label{propproj}
Under~(\ref{IRMproj}), if $W$ and $\hat\Sigma$ are regular,
then the maximum likelihood estimator of $(\mu,V,b,c)$ is defined by:
\begin{itemize}
\item $\hat b$ is the eigenvector
associated to the largest eigenvalue $\hat\lambda$ 
of $\hat \Sigma^{-1} M^t W^{-1} M$,
\item $\hat c = \displaystyle \frac{1}{\hat b^t \hat V \hat b} W^{-1} M \hat b$,
\item $\hat \mu = \bar{X} - \bar{s}^t \hat c \hat V \hat b$,
\item  $\hat V= \hat\Sigma -\displaystyle\frac{\hat\lambda}{{\hat b}^t\hat\Sigma \hat b}
\hat\Sigma \hat b {\hat b}^t \hat \Sigma$.
\end{itemize}
\end{proposition}
As a consequence of the above equation, one also has
$\hat \lambda= 1 - {\hat b}^t \hat V \hat b{/{\hat b}^t \hat \Sigma \hat b}$,
which provides an estimation of the signal to noise ratio
in the direction $b$ through 
\begin{equation}
\label{lambdahat}
\hat\rho=\hat\lambda/(1-\hat\lambda).
\end{equation}
Let us now show that the SIR method corresponds
to the particular case of piecewise constant basis functions
$s_j(.)$, $j=1,\dots,h$.

\subsection{Sliced Inverse Regression (SIR)}
%%%%%%%%%%%%%%%%%%%%%%%%%%%%%%%%%%%%%%%%%%%%%%%%%%%%%%%%%%%%%%%%%%%%%%%
\label{SIR}

Suppose the range of $Y$ is partitioned into $h+1$ non-overlapping
slices $S_j$, $j=1,\dots,h+1$ and consider the $h$ basis functions defined by
\begin{equation}
\label{baseindic}
s_j(.)=\I\{.\in S_j\},\;\; j=1,\dots,h
\end{equation}
where $\I$ is the indicator function.
Let us denote by $n_j$ the number of $Y_i$ in slice 
$j=1,\dots,h+1$, 
define the corresponding proportion by $f_j=n_j/n$,
the empirical mean of $X$ given $Y\in S_j$ by
$$
\bar{X}_j = \frac{1}{n_j} \sum_{Y_i\in S_j} X_i
$$
and let $\hat\Gamma$ be the $p\times p$ empirical "between slices"
covariance matrix defined by
$$
\hat\Gamma=\sum_{j=1}^{h+1} f_j (\bar{X}_j-\bar{X})  (\bar{X}_j-\bar{X})^t.
$$
In this context, the following consequence of Proposition~\ref{propproj}
can be established.
\begin{corollary}
\label{corodiscr}
Under~(\ref{IRMproj}) and~(\ref{baseindic}), if $\hat\Sigma$ is regular,
then the maximum likelihood estimator $\hat b$ of $b$ is the eigenvector
associated to the largest eigenvalue of ${\hat\Sigma}^{-1} \hat\Gamma$.
\end{corollary}
To summarize, the maximum likelihood estimator of $b$ is 
(up to a scale factor) the 
classical SIR estimator of the direction $\beta$ spanning the central
subspace. 
The next section is dedicated
to the introduction of a regularization in the inverse regression
problem in order to avoid the inversion of $\hat\Sigma$.

%%%%%%%%%%%%%%%%%%%%%%%%%%%%%%%%%%%%%%%%%%%%%%%%%%%%%%%%%%%%%%%%%%%%%%%
\section{Regularized inverse regression}
%%%%%%%%%%%%%%%%%%%%%%%%%%%%%%%%%%%%%%%%%%%%%%%%%%%%%%%%%%%%%%%%%%%%%%%
\label{RIR}

First, we present in Subsection~\ref{oa} how the introduction of
a prior information on the unknown direction $b$ can overcome the SIR
limitations due to the ill-conditioning or singularity of $\hat\Sigma$.
Second, some links with the existing SIR
regularizations are highlighted in Subsection~\ref{links}.
Finally, basing on our framework, three new regularizations of the
SIR method are introduced in Subsection~\ref{new}.

\subsection{Introducing a Gaussian prior}
%%%%%%%%%%%%%%%%%%%%%%%%%%%%%%%%%%%%%%%%%%%%%%%%%%%%%%%%%%%%%%%%%%%%%%%
\label{oa}

We propose to introduce a prior information on the projection of $X$ on $b$ appearing in the inverse regression model. 
More precisely, we focus on 
$$
\frac{\E(b^t (X-\mu)|Y)b}{b^t V b} =  c(Y)  b =  s^t(Y) c b ,
$$
from~(\ref{eqc}) and (\ref{projc}), assuming that
\begin{equation}
\label{apriori}
\left( 1+\rho \right)^{-1/2}
(s(Y)-\bar{s})^t c b\sim{\cal N}(0,\Omega).
\end{equation}
The role of the matrix $\Omega$ is to describe which directions in
$\R^p$ are the most likely to contain $b$. Some examples are
provided in the next two paragraphs.
The scalar $(1+\rho)^{-1/2}$ is introduced for normalization purposes,
permitting to preserve the interpretation of the eigenvalue in
terms of signal to noise ratio.
As a comparison, in~\cite{Anes2}, a Bayesian estimation method
is proposed basing on B-splines approximation of the link function $g$
in the forward model and a Fisher-von Mises prior on the direction $b$. 
In our approach, working on the inverse regression model
allows to obtain explicit solutions, see Lemma~\ref{lemmap} and 
Proposition~\ref{propregul} below.
\begin{lemma}
\label{lemmap}
Maximum likelihood estimators are obtained by minimizing
\begin{equation}
\label{Gregul}
G_\Omega(\mu,V,b,c)= G(\mu,V,b,c)
+ \frac{(b^t \Omega^{-1} b)(b^t V b)(c^t W c) }{b^t\Sigma b}
\end{equation}
with respect to $(\mu,V,b,c)$. 
\end{lemma}
Comparing to~(\ref{Gproj}), the additional
term due to the prior information can be read as a 
regularization term in Tikhonov theory, see for instance~\cite{Vogel},
Chapter~1, penalizing large projections. 
The following result can be stated:
\begin{proposition}
\label{propregul}
Under~(\ref{IRMproj}) and~(\ref{apriori}), if $W$
and $\Omega\hat\Sigma+I_p$
are regular, then the maximum likelihood estimator of 
$(\mu,V,b,c)$ is defined by:
\begin{itemize}
\item $\hat b$ is the eigenvector
associated to the largest eigenvalue $\hat\lambda$ 
of  $(\Omega\hat\Sigma+I_p)^{-1} \Omega M^t W^{-1} M$,
\item $\hat c = \displaystyle \frac{1}{(1+\eta(\hat b))\hat b^t \hat V \hat b}
W^{-1} M \hat b$,
where $\eta(\hat b)=\displaystyle \frac{{\hat b}^t \Omega^{-1} \hat b}{ {\hat b}^t \hat \Sigma  \hat b}$,
\item $\hat \mu = \bar{X} - \bar{s}^t \hat c \hat V \hat b$,
\item $\hat V= \hat\Sigma -\displaystyle\frac{\hat\lambda}{{\hat b}^t\hat\Sigma \hat b}
\hat\Sigma \hat b {\hat b}^t \hat \Sigma$.
\end{itemize}
\end{proposition}
Compared to Proposition~\ref{propproj}, 
the inversion of $\hat\Sigma$ is replaced by the inversion
of $\Omega\hat\Sigma+I_p$. Thus, for a properly chosen prior matrix
$\Omega$, the numerical instabilities in the estimation of $b$ disappear.
Note that, since the estimation of $V$ is formally the same as in 
Proposition~\ref{propproj}, the interpretation (\ref{lambdahat}) of 
$\hat \lambda$ still holds.
As previously, this result can be applied to the particular case
of the SIR method.

\begin{corollary}
\label{cororegul}
Under~(\ref{IRMproj}), (\ref{baseindic}) and~(\ref{apriori}), 
if $\Omega\hat\Sigma+I_p$ is regular,
then the maximum likelihood estimator of $b$ 
is the eigenvector
associated to the largest eigenvalue of $(\Omega\hat\Sigma+ I_p)^{-1} \Omega\hat\Gamma$.
\end{corollary}
In the following, the above estimator of the direction $b$ will be referred to as
the Gaussian Regularized Sliced Inverse Regression (GRSIR) estimator.
Let us emphasize that the GRSIR estimator can be extended to the
multi-index situation by considering the eigenvectors $b_1,\dots,b_K$
associated to the $K$ largest eigenvalues $\lambda_1>\dots>\lambda_K$
of $(\Omega\hat\Sigma+ I_p)^{-1} \Omega\hat\Gamma$.
For instance, one can show that $b_2$ maximizes $G_\Omega$ under 
the orthogonality
constraint $b_1^t (\Sigma+\Omega^{-1}) b_2=0$.

Some examples of possible prior covariance matrices $\Omega$ are now
presented.

\subsection{Links with existing methods}
%%%%%%%%%%%%%%%%%%%%%%%%%%%%%%%%%%%%%%%%%%%%%%%%%%%%%%%%%%%%%%%%%%%%%%%
\label{links}

In all the next examples, a non-negative parameter $\tau$
is introduced in the prior covariance matrix in order to 
tune the importance of the penalty term in~(\ref{Gregul}). 
Consequently, in the sequel, $\tau$ is called a regularization
parameter.

\paragraph{Classical SIR approach.}
It is easily seen from Corollary~\ref{cororegul} that
choosing the prior covariance matrix
$\Omega_0=(\tau\hat\Sigma)^{-1}$ in GRSIR gives back 
SIR, and this for all $\tau>0$.
This prior matrix indicates that directions
corresponding to small variances are most likely, {\it i.e}
the SIR method favors directions in which $\hat\Sigma$
is close to singularity. In practice, this choice yields
instabilities in the estimation.

\paragraph{Ridge approach~\cite{Scrucca,Scrucca2,Zhong}.}
The simplest choice for the prior covariance matrix is $\Omega_1=\tau^{-1} I_p$.
In this case, the identity matrix indicates that no privileged direction
for $b$ is available. 
Following Corollary~\ref{cororegul},
the corresponding GRSIR estimator of $b$ is the eigenvector of $(\hat\Sigma+\tau I_p)^{-1} \hat\Gamma$
associated to its largest eigenvalue, which is the ridge estimator introduced
independently in~\cite{Zhong} and~\cite{Scrucca,Scrucca2}. 

\paragraph{PCA+SIR approach~\cite{Chiar,LiLi}.}
As already seen, a popular technique to overcome the singularity problems
of $\hat\Sigma$
is to use PCA as a preprocessing step~\cite{Chiar,LiLi}. The principle
is the following. Let $d\in\{1,\dots,p\}$ be fixed and denote by $\lambda_1\geq
\dots\geq \lambda_d$ the $d$ largest eigenvalues of $\hat\Sigma$
(supposed to be positive), $q_1,\dots,q_d$ the associated
eigenvectors and $S_d=\mbox{span}(q_1,\dots,q_d)$
the linear subspace spanned by $q_1,\dots,q_d$.
The first step consists in projecting the predictors on $S_d$.
The second step is to perform SIR in this subspace. 
The next result shows that this method corresponds to a particular
prior covariance matrix.
\begin{proposition}
\label{propacpsir}
PCA+SIR corresponds to GRSIR with prior covariance matrix
$$
\Omega_2=\frac{1}{\tau}\sum_{j=1}^d \frac{1}{\lambda_j} q_j q_j^t,
$$
where $\tau>0$ is arbitrary.
\end{proposition}
Let us note that although $\Omega_2$ depends on $\tau$,
the GRSIR estimator does not, since there is no regularization parameter
in the PCA+SIR methodology.

\paragraph{Li and Yin's approach~\cite{LiYin}.} Their proposed 
$L_2-$ regularization consists in estimating $(b,c)$ by minimization of
$$
H_\tau(b,c)=\sum_{j=1}^h f_j (\mu+\hat\Sigma c_j b-\bar{X}_j )^t N
(\mu+\hat\Sigma c_j b-\bar{X}_j )+ \tau b^t b,
$$
the matrix $N$ being either $N=I_p$ or $N=\hat\Sigma^{-1}$.
In our opinion, this approach suffers from a lack of invariance since
the functional $H_\tau$ does not penalize the same way two different
axes ($b$ and $2b$ for instance) defining the same direction.
As a consequence, we have shown (\cite{Biometrics}, Proposition~1)
that the only possible solution $\hat b$
of the minimization problem is $\hat b=0$:
In view of this result, the proposed alternating least squares
algorithm (\cite{LiYin}, Section~2) cannot be justified theoretically.
As a comparison, 
our method does not yield this 
kind of problem thanks to the invariance property:
$G_\Omega(\mu,V,tb,c/t)=G_\Omega(\mu,V,b,c)$
for all real numbers $t\neq 0$.
We now propose some alternative choices of the covariance 
matrix $\Omega$ yielding new regularizations of the SIR method.

\subsection{Three new SIR regularizations}
%%%%%%%%%%%%%%%%%%%%%%%%%%%%%%%%%%%%%%%%%%%%%%%%%%%%%%%%%%%%%%%%%%%%%%%
\label{new}

\paragraph{Tikhonov regularization.}
An alternative choice of the prior covariance matrix is 
$\Omega_3=\tau^{-1} \hat\Sigma$. Comparing $\Omega_3$ to the matrix $\Omega_0$
associated to the SIR method, it appears that the underlying ideas
are opposite. Here, directions
corresponding to large variances are most likely. 
The associated GRSIR estimator of the direction $b$
is the eigenvector of $(\hat\Sigma^2+\tau I_p)^{-1} \hat\Sigma\hat\Gamma$
associated to its largest eigenvalue. In the following,
this estimator will be referred to as the Tikhonov estimator.
Indeed, let us recall that the classical SIR estimator 
is obtained by a spectral decomposition
of $\hat\Sigma^{-1}\hat\Gamma$. For all $k=1,\dots,p$, denote by $x_k$
the $k-$th column of this matrix.
Computing $x_k$ is equivalent to solving with respect to $x$ the linear system
$\hat\Sigma x = \hat\Gamma_k$ where $\hat\Gamma_k$
is the $k-$th column of $\hat\Gamma$. 
The associated Tikhonov minimization problem (see (1.34) in~\cite{Vogel}) 
can be written as
$$
x_k=\arg\min_x \| \hat\Sigma x - \hat\Gamma_k\|^2 + \tau \|x\|^2
= (\hat\Sigma^2+\tau I_p)^{-1}\hat\Sigma\hat\Gamma_k.
$$
Thus, in this framework, $(\hat\Sigma^{-1}\hat\Gamma)_k$ is estimated
by $(\hat\Sigma^2+\tau I_p)^{-1}\hat\Sigma\hat\Gamma_k$ and consequently
$\hat\Sigma^{-1}\hat\Gamma$ is estimated by $(\hat\Sigma^2+\tau I_p)^{-1}\hat\Sigma\hat\Gamma$.

\paragraph{Dimension reduction approaches.}
It has been seen in Proposition~\ref{propacpsir}, that the PCA+SIR
approach is equivalent to using the prior covariance matrix $\Omega_2$ in 
GRSIR.  The following result is an extension
to more general covariance matrices.
\begin{proposition}
\label{propreducdim}
For all real function $\varphi$ let
$$
\Omega(\varphi)=\sum_{j=1}^d \varphi(\lambda_j) q_j q_j^t.
$$
Then, the associated GRSIR estimator
can be obtained by first projecting the predictors on  $S_d=\mbox{span}(q_1,\dots,q_d)$ and
second performing GRSIR on the projected
predictors with prior covariance matrix
$$
\tilde\Omega(\varphi)=\sum_{j=1}^p \varphi(\lambda_j) q_j q_j^t.
$$
\end{proposition}
The dimension $d$ plays the role of a "cut-off" parameter, since
when computing $\hat b$, all directions $q_{d+1},\dots,q_p$ are
discarded.
Three illustrations of this result can be given:
\begin{itemize}
\item Choosing $\varphi(t)=1/(\tau t)$, we obtain
$\Omega(1/(\tau \Id))=\Omega_2$ and 
$\tilde\Omega(1/(\tau \Id))=(\tau\hat\Sigma)^{-1}=\Omega_0$,
where $\Id$ is the identity function.
It appears that Proposition~\ref{propacpsir} is 
a particular case of Proposition~\ref{propreducdim}.
As already discussed, since the choice of $\Omega_0$ as a 
prior covariance matrix seems not very natural,
we thus propose two new choices. 
\item First, $\varphi(t)=1/\tau$ yields 
$$
\Omega_4\stackrel{def}{=}\Omega(1/\tau)=\frac{1}{\tau}\sum_{j=1}^d q_j q_j^t,
$$
and 
$\tilde\Omega(1/\tau)=I_p/\tau=\Omega_1$. 
Consequently, this new method consists in applying 
the ridge  approach~\cite{Zhong} on the projected predictors,
the interpretation being that, in the subspace $S_d$, all directions
share the same prior probability. This method will be referred to 
as PCA+ridge.

\item Second, $\varphi(t)=t/\tau$ yields 
$$
\Omega_5\stackrel{def}{=}\Omega(\Id/\tau)=\frac{1}{\tau}\sum_{j=1}^d \lambda_j
q_j q_j^t,
$$
and 
$\tilde\Omega(\Id/\tau)=\hat\Sigma/\tau=\Omega_3$. This new method consists in applying 
Tikhonov approach on the projected predictors. In this context,
directions of $S_d$ carrying a large fraction of the total variance 
of $X$ are more likely. This method will be referred to as PCA+Tikhonov.
\end{itemize}

%%%%%%%%%%%%%%%%%%%%%%%%%%%%%%%%%%%%%%%%%%%%%%%%%%%%%%%%%%%%%%%%%%%%%%%
\section{Numerical experiments}
%%%%%%%%%%%%%%%%%%%%%%%%%%%%%%%%%%%%%%%%%%%%%%%%%%%%%%%%%%%%%%%%%%%%%%%
\label{FSP}

GRSIR methods associated to the 
prior covariance matrices $\Omega_0$ (SIR), $\Omega_1$ (ridge),
$\Omega_2$ (PCA+SIR), $\Omega_3$ (Tikhonov), $\Omega_4$ (PCA+ridge)
and $\Omega_5$ (PCA+Tikhonov) are compared on simulated data.
A sample of size $n=100$ of the random 
pair $(X,Y)$ is considered, where $X \in\R^p$ with 
$p=50$ and $Y\in \R$.
The random vector $X$ is Gaussian, centered, with covariance matrix
$\Sigma=Q \Delta Q^t$ where 
$\Delta$ is the diagonal matrix
containing the eigenvalues of $\Sigma$ defined by
$\Delta=$diag$\{p^\theta,(p-1)^\theta,\dots,1^\theta\}$ and
$Q$ is a randomly chosen orthogonal
matrix. Several values of $\theta$ will be considered. Note
that the condition number of $\Sigma$ is given by $p^\theta$
and is thus an increasing function of $\theta$. As suggested by~\cite{Heib}, the orthogonal matrix $Q$ is obtained as follow: First, we construct a $p \times p$ matrix where each coefficients are randomly chosen from a Gaussian distribution. Next, the matrix $Q$ is obtained by applying the QR-decomposition on this matrix. We consider two SIR models: 
$$
{\bf{Model \ 1 \ \ }} \ Y=\sin\left(\frac{\pi}{2\sigma} \beta^t X\right)+\varepsilon,
$$
$$
{\bf{Model \ 2 \ \ }} \ Y= \left | \frac{\beta^t X}{\sigma} -1/2 \right |+\varepsilon,
$$
where $\sigma$ is the standard deviation of the projection of $X$
on $\beta$ {\it i.e.} $\sigma=(\beta^t\Sigma\beta)^{1/2}$, $\varepsilon$ is a centered Gaussian random value with standard deviation 0.03 independent of $X$.  The true index is $\beta=5^{-1/2}Q(1,1,1,1,1,0,\dots,0)^t$. Hence in the direction of $\beta$, the variance of $X$ is large.

In all the experiments described below, we replicate $N=100$ times Models~1 and~2. More precisely, for Model 1, we compute the pairs $(X,Y^{(r)})$, $r=1,\ldots,N$ where
$$
Y^{(r)} = \sin\left(\frac{\pi}{2\sigma} \beta^t X\right)+\varepsilon_r, \ \varepsilon_r \sim {\cal{N}}(0,(0.03)^2),  
$$
with $\varepsilon_r$ independent of $X$. The same is done for Model 2. We then obtained $N$ estimators ${\hat{b}}^{(r)}$, $r=1,\ldots,N$ of $\beta$. In order to evaluate the quality of the estimate ${\hat{b}}$, we compute two criteria. The first one is the mean of the squared cosine (MSC):
$$
{\mbox{MSC}} = \frac{1}{N} \sum_{r=1}^N (\beta^t{\hat{b}}^{(r)})^2. 
$$
The second one is a kind of variance of the squared cosine (VSC):
$$
{\mbox{VSC}} = \frac{1}{N(N-1)} \sum_{r=1}^{N} \sum_{r \neq s} (({\hat{b}}^{(s)})^t{\hat{b}}^{(r)})^2.
$$
These criteria can be adapted to multi-index models~\cite{Ferre}.
The closer these quantities are to 1, the better the estimation is.

\begin{figure}[h]
\begin{center}
\vspace*{1cm}
\begin{minipage}{0.4\textwidth}
\begin{center}
\epsfig{figure=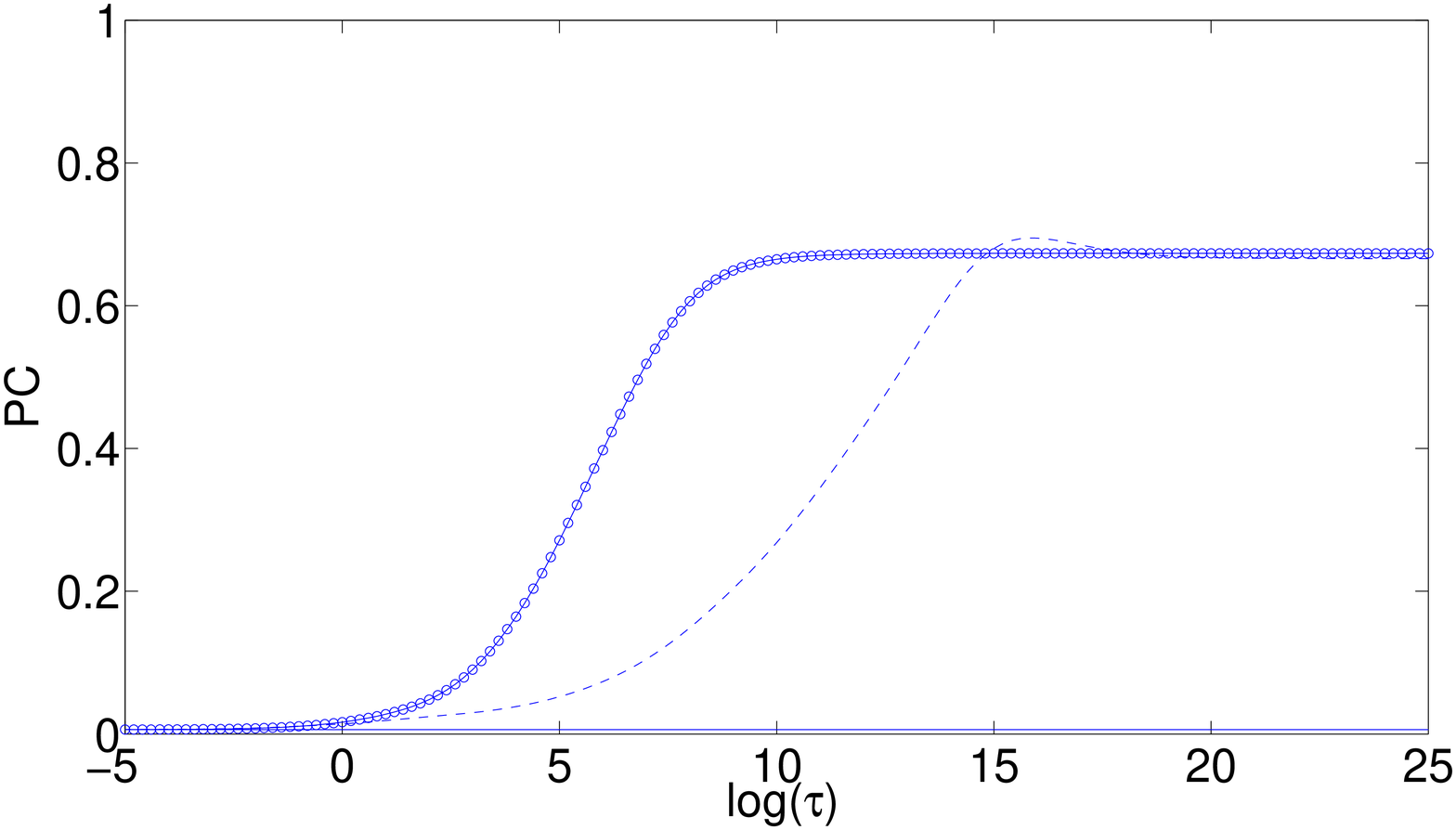,width=1\textwidth,height=1\textwidth}
\centerline{(a) Criterion MSC}
\end{center}
\end{minipage}
\begin{minipage}{0.4\textwidth}
\epsfig{figure=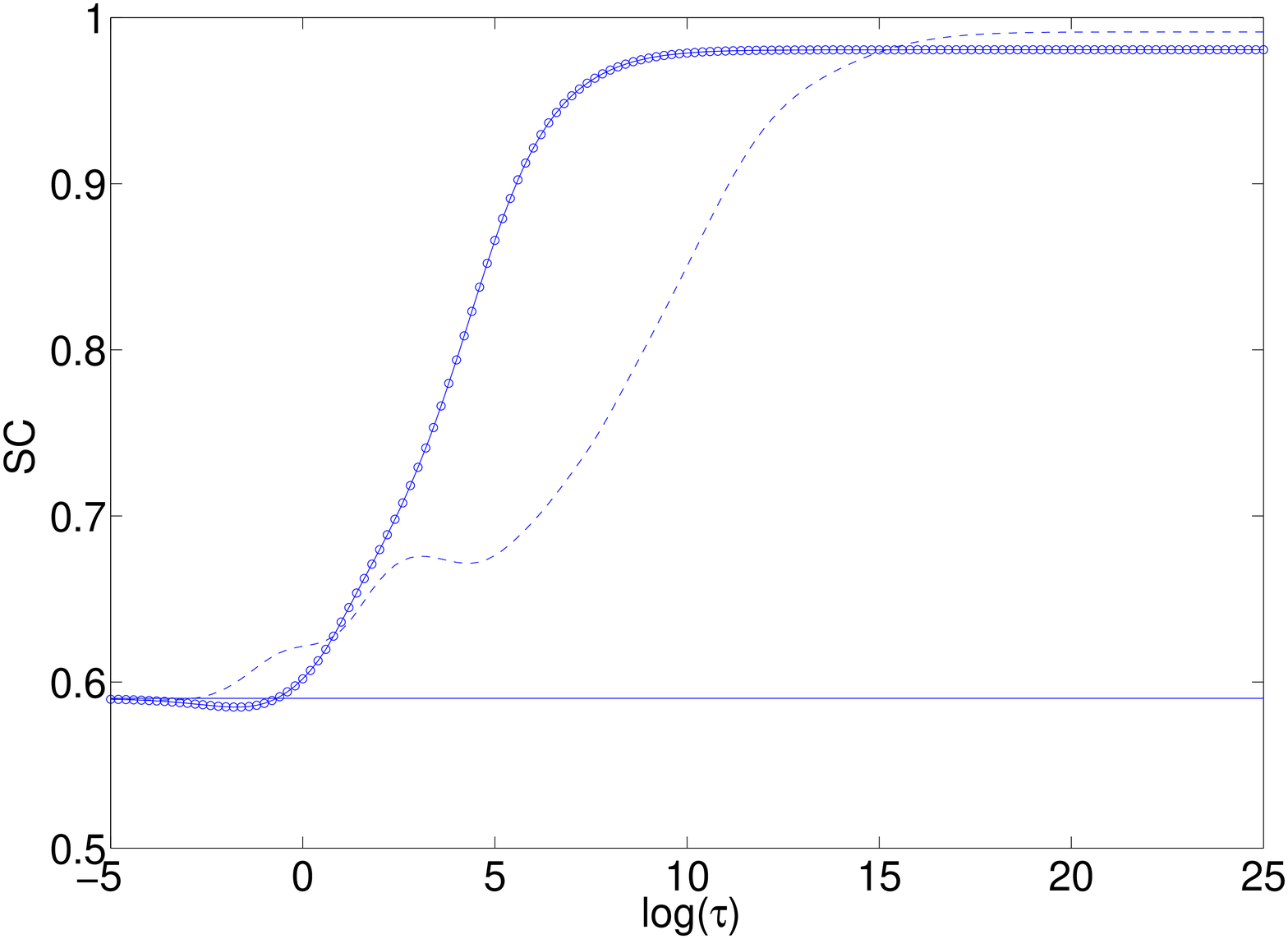,width=1\textwidth,height=1\textwidth}
\centerline{(b) Criterion VSC}
\end{minipage}
\vspace*{1cm}
\caption{Values of MSC and SSC with respect to the regularization parameter for Model 1. The condition number is fixed to $\theta=2$.
Horizontally: $\log(\tau)$, vertically: (a) MSC and (b) VSC.
Continuous line: $\Omega_0$ (SIR), "-o-" line: $\Omega_1$ (ridge),
dashed line: $\Omega_3$ (Tikhonov).}
\label{figzhotik}
\end{center}
\end{figure}

\begin{figure}[h]
\begin{center}
\vspace*{1cm}
\begin{minipage}{0.4\textwidth}
\epsfig{figure=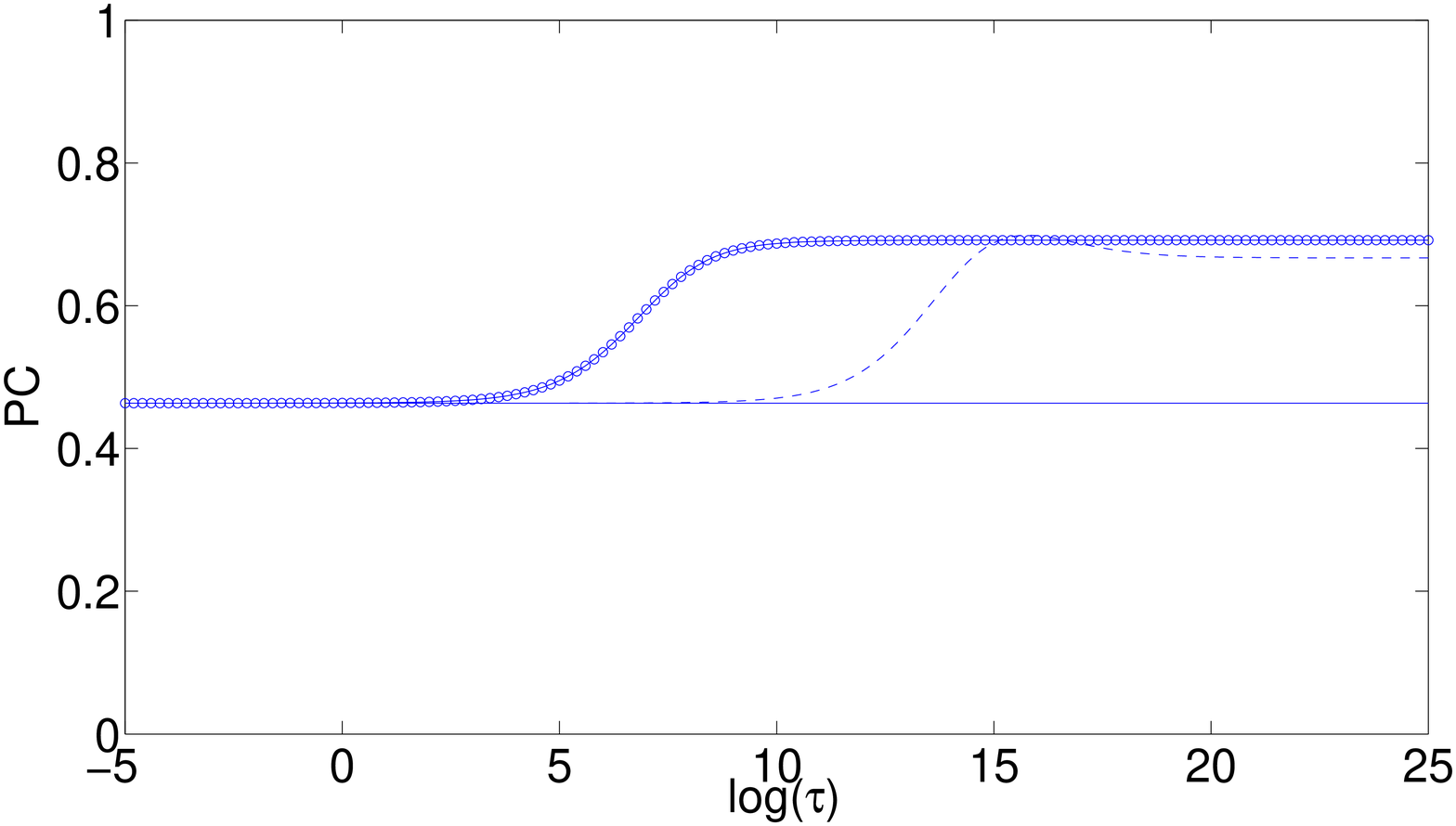,width=1\textwidth,height=1\textwidth}
\centerline{(a) Criterion MSC}
\end{minipage}
\begin{minipage}{0.4\textwidth}
\epsfig{figure=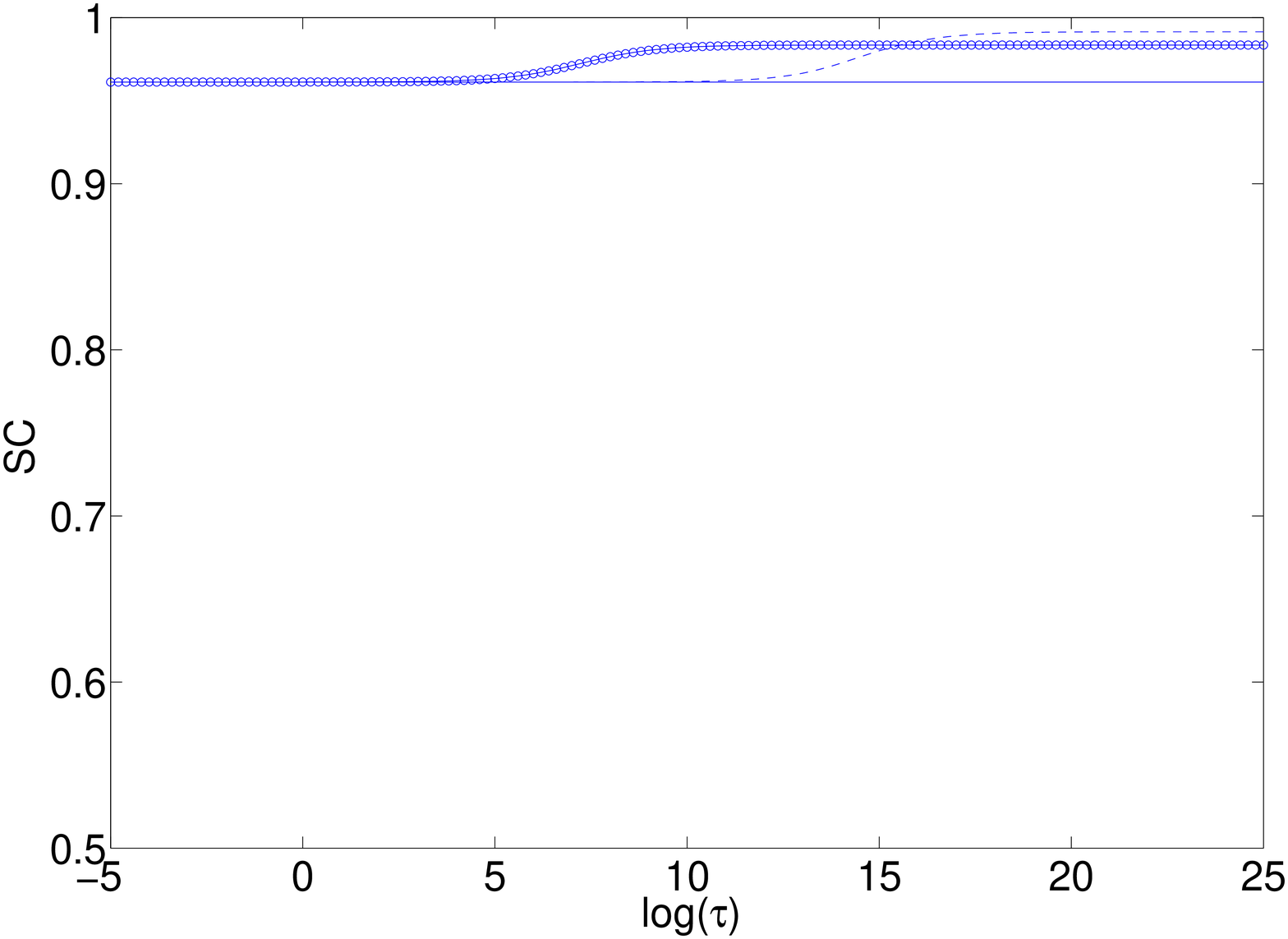,width=1\textwidth,height=1\textwidth}
\centerline{(b) Criterion VSC}
\end{minipage}
\vspace*{1cm}
\caption{Values of MSC and VSC with respect to the regularization parameter for Model 1. The cut-off dimension is chosen to $d=20$ and the condition number is fixed to $\theta=2$.
Horizontally: $\log(\tau)$, vertically: (a) MSC and (b) VSC.
Continuous line: $\Omega_2$ (PCA+SIR), "-o-" line: $\Omega_4$ (PCA+ridge),
dashed line: $\Omega_5$ (PCA+Tikhonov).}
\label{figacpzhotik}
\end{center}
\end{figure}

\begin{figure}[h]
\begin{center}
\vspace*{1cm}
\begin{minipage}{0.4\textwidth}
\begin{center}
\epsfig{figure=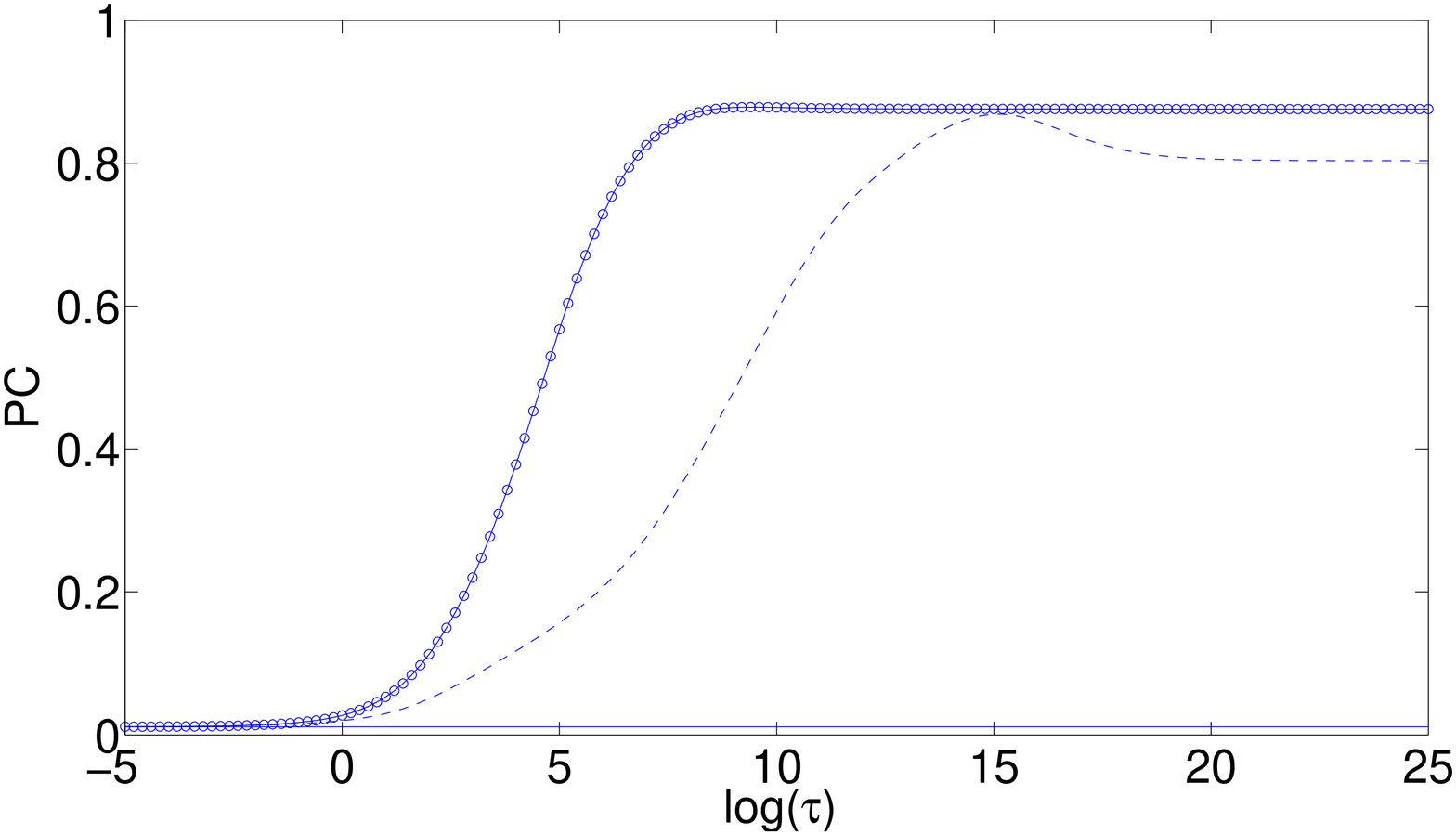,width=1\textwidth,height=1\textwidth}
\centerline{(a) Criterion MSC}
\end{center}
\end{minipage}
\begin{minipage}{0.4\textwidth}
\epsfig{figure=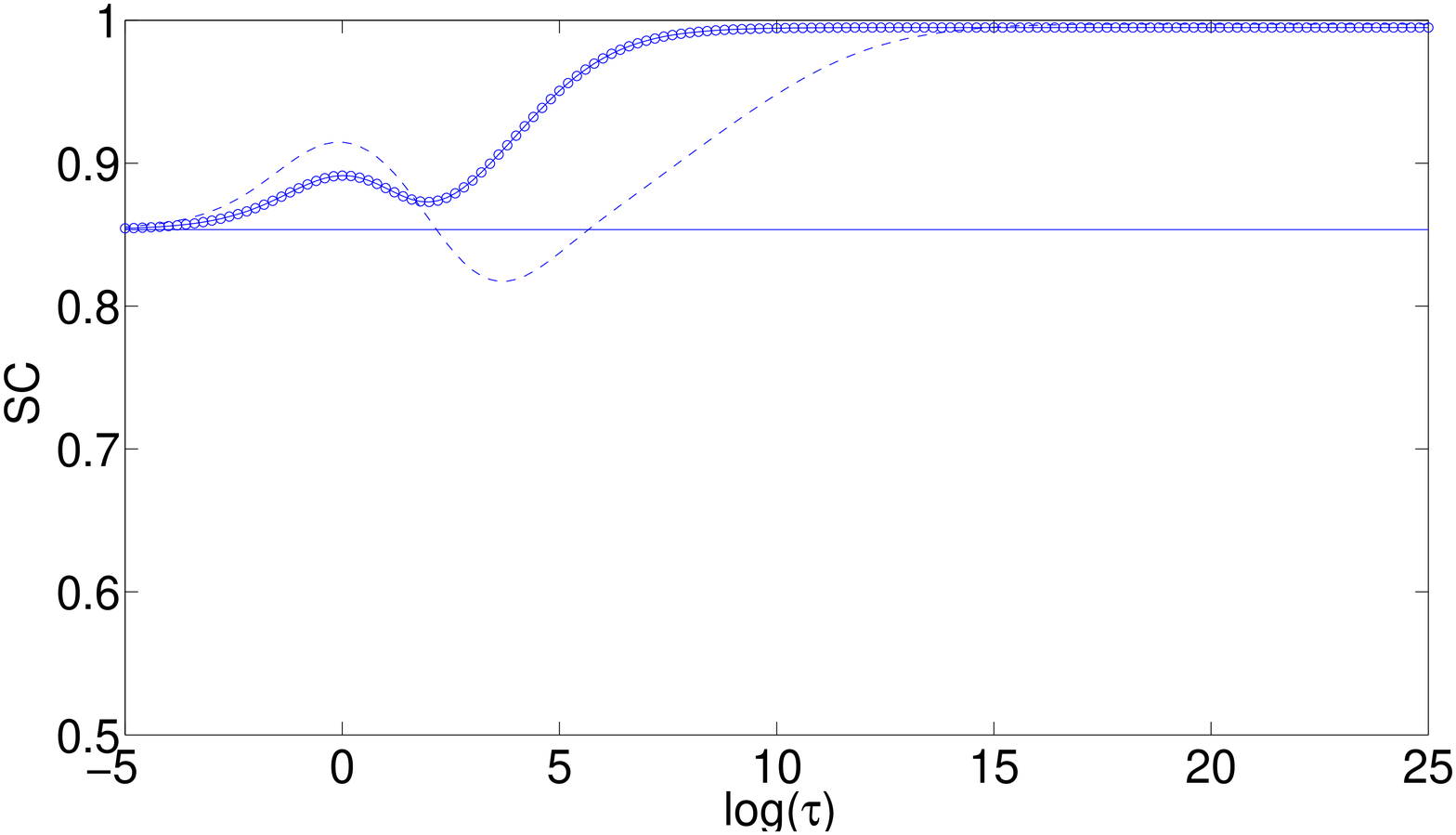,width=1\textwidth,height=1\textwidth}
\centerline{(b) Criterion VSC}
\end{minipage}
\vspace*{1cm}
\caption{Values of MSC and VSC with respect to the regularization parameter for Model 2. The condition number is fixed to $\theta=2$.
Horizontally: $\log(\tau)$, vertically: (a) MSC and (b) VSC.
Continuous line: $\Omega_0$ (SIR), "-o-" line: $\Omega_1$ (ridge),
dashed line: $\Omega_3$ (Tikhonov).}
\label{figzhotiknew}
\end{center}
\end{figure}

\begin{figure}[h]
\begin{center}
\vspace*{1cm}
\begin{minipage}{0.4\textwidth}
\epsfig{figure=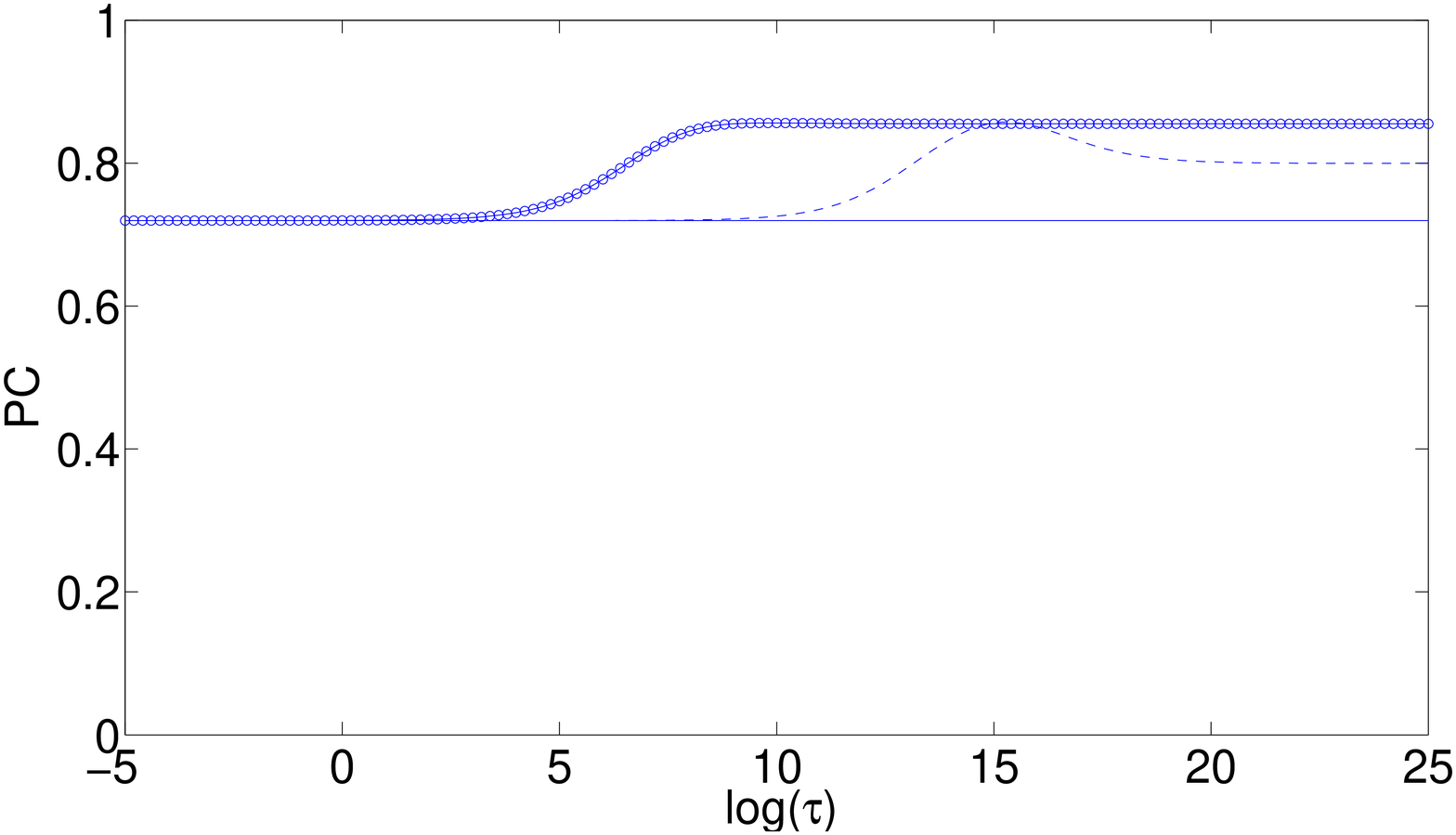,width=1\textwidth,height=1\textwidth}
\centerline{(a) Criterion MSC}
\end{minipage}
\begin{minipage}{0.4\textwidth}
\epsfig{figure=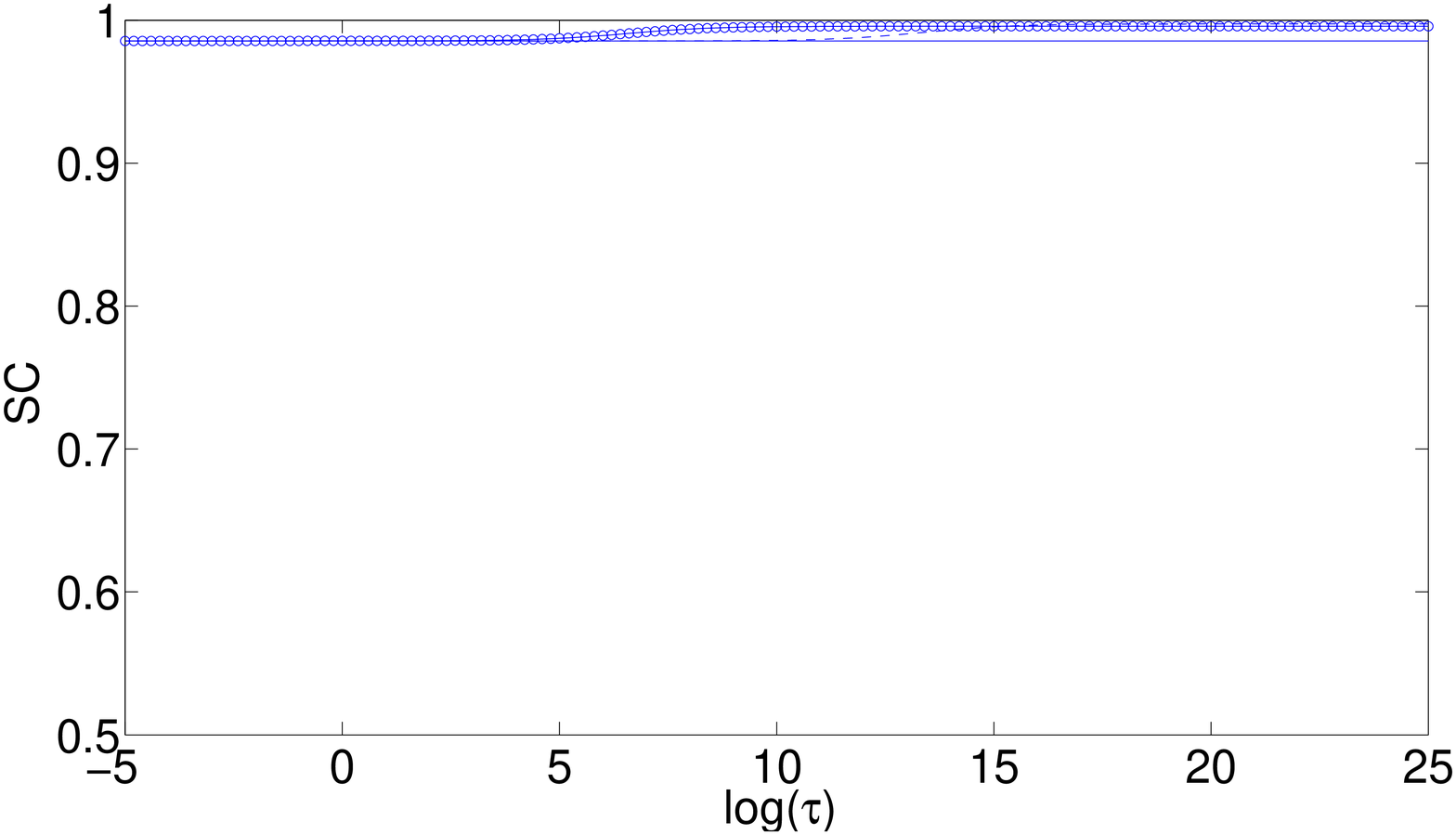,width=1\textwidth,height=1\textwidth}
\centerline{(b) Criterion VSC}
\end{minipage}
\vspace*{1cm}
\caption{Values of MSC and VSC with respect to the regularization parameter for Model 2. The cut-off dimension is chosen to $d=20$ and the condition number is fixed to $\theta=2$.
Horizontally: $\log(\tau)$, vertically: (a) MSC and (b) VSC.
Continuous line: $\Omega_2$ (PCA+SIR), "-o-" line: $\Omega_4$ (PCA+ridge),
dashed line: $\Omega_5$ (PCA+Tikhonov).}
\label{figacpzhotiknew}
\end{center}
\end{figure}

\begin{figure}[p]
\begin{center}
\epsfig{figure=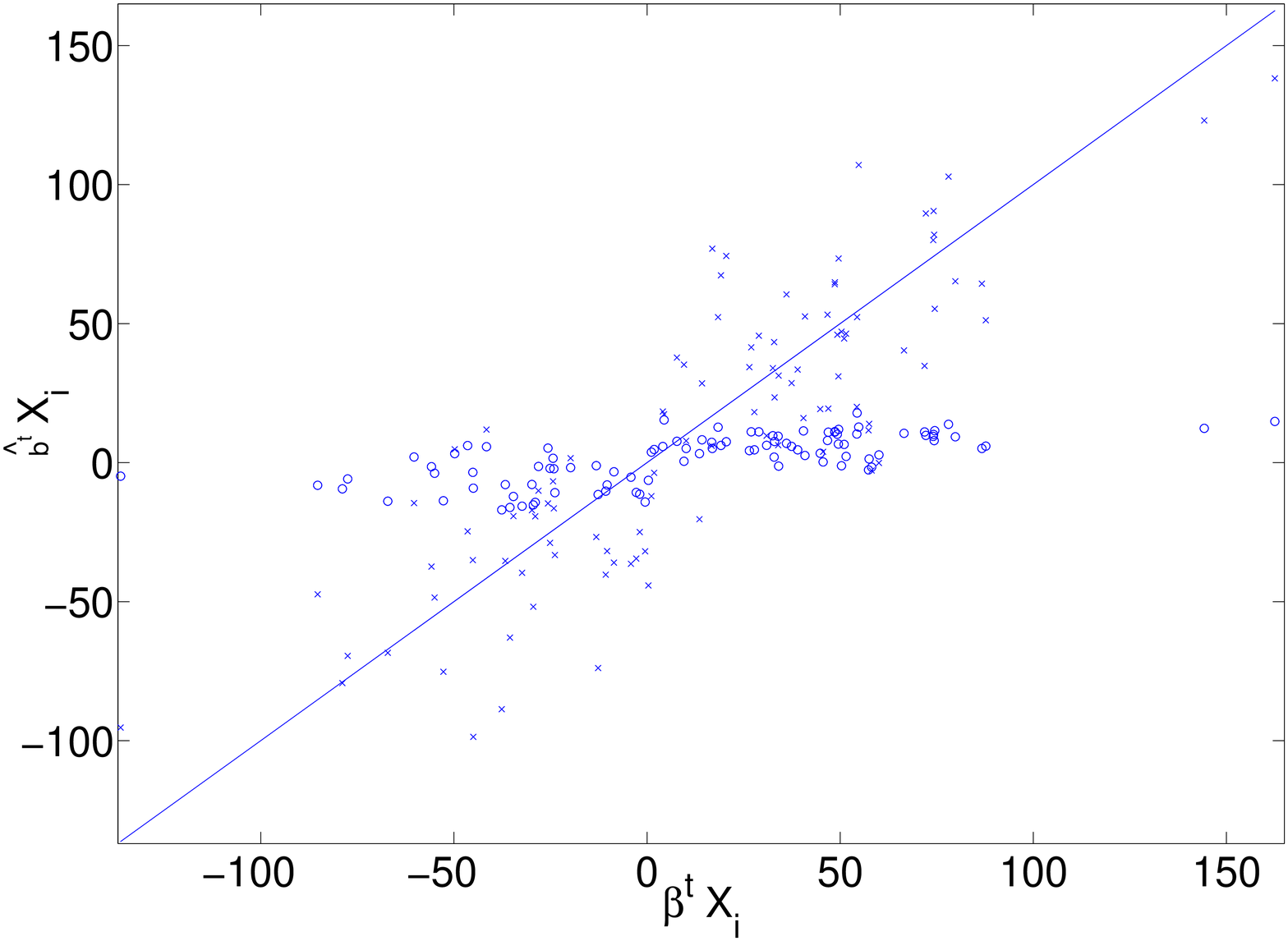,width=0.6\textwidth,height=0.5\textwidth}
\end{center}
\caption{An example of pairs $(\beta^t X_i, {\hat{b}}^t X_i)$, $i=1,\dots,n$ for Model 1 where the axis ${\hat{b}}$ is computed with the classical SIR method (points "o") and the PCA+Tikhonov method (points "x"). Here, $\theta=2$ and $d=20$.}
\label{compasirtik}
\end{figure}

\begin{figure}[p]
\begin{center}
\epsfig{figure=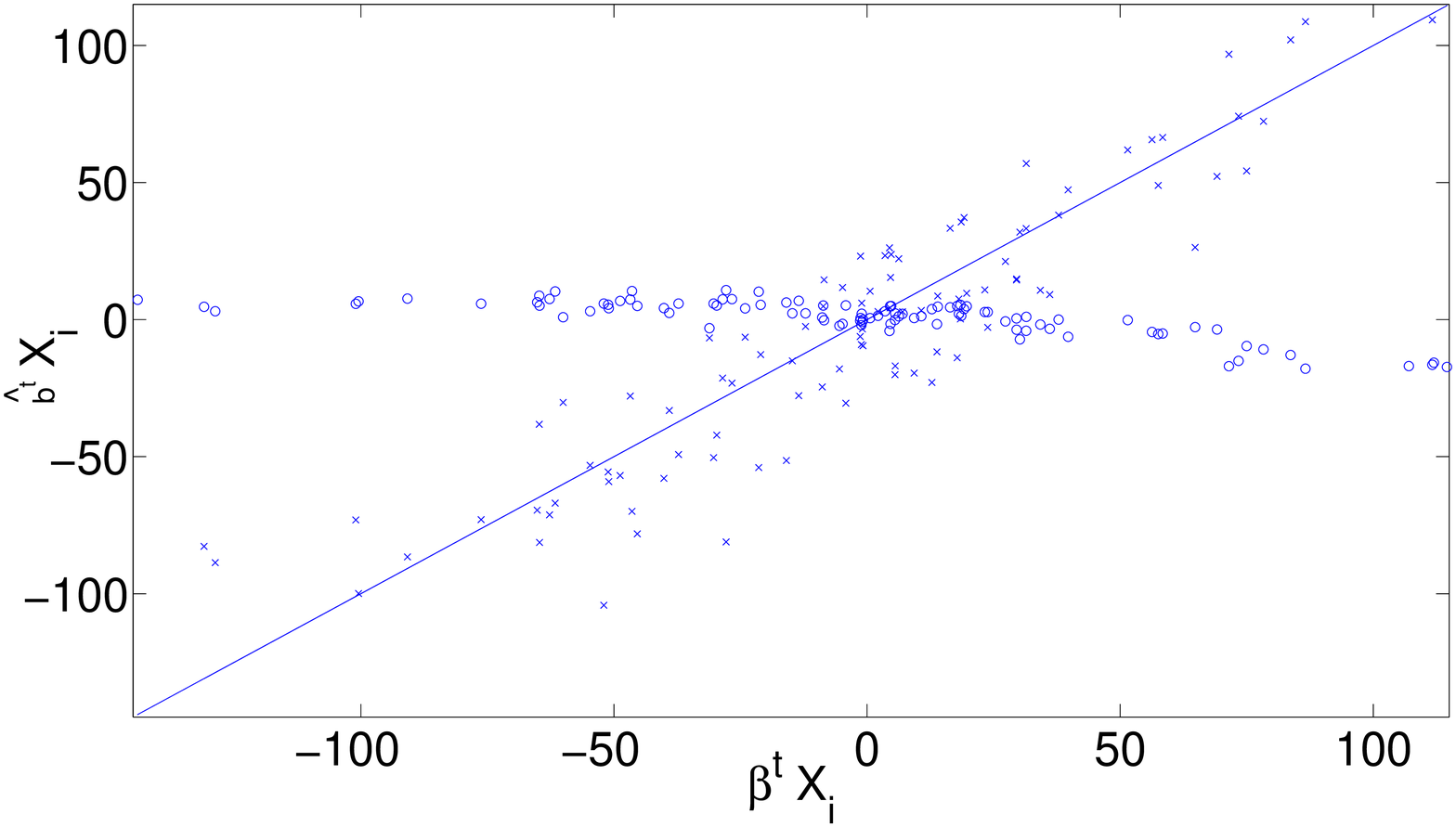,width=0.6\textwidth,height=0.5\textwidth}
\end{center}
\caption{An example of pairs $(\beta^t X_i, {\hat{b}}^t X_i)$, $i=1,\dots,n$ for Model 2 where the axis ${\hat{b}}$ is computed with the classical SIR method (points "o") and the PCA+Tikhonov method (points "x"). Here, $\theta=2$ and $d=20$.}
\label{compasirtiknew}
\end{figure}

First experiment: We illustrate the dependence of the above methods with respect to the regularization
parameter $\tau$. In this aim, we compute the two previous defined criteria as function of the regularization parameter. A logarithmic scale was adopted, 
150 values of $\log(\tau)$ regularly distributed in $[-5,25]$
were considered.
Here, we limit ourselves to $\theta=2$. Moreover, we choose
$d=20$ in the PCA+SIR, PCA+ridge and PCA+Tikhonov methods.
Results are displayed on Figure~\ref{figzhotik} and Figure~\ref{figacpzhotik} for Model~1 and on Figure~\ref{figzhotiknew} and Figure~\ref{figacpzhotiknew} for Model~2.
It appears that the classical SIR approach gives very poor results
in such a situation where $n$ and $p$ are of the same order ($n/p=2$).
Ridge and Tikhonov regularizations can bring a significant improvement
provided $\tau$ is large enough. The VSC criterion of the ridge and the Tikhonov methods is better than the one of the classical SIR method if the regularization parameter is taken sufficiently large.
PCA+SIR obtains reasonable results compared to SIR, with the
advantage of do not requiring the selection of $\tau$.
The selection of $d$ is addressed in our third experiment.
Note that PCA+ridge and PCA+Tikhonov methods are less sensitive
to the choice of $\tau$ than ridge and Tikhonov methods. Concerning the criterion VSC, PCA+ridge and PCA+Tikhonov methods both outperform the PCA+SIR for sufficiently large values of the regularization parameter.
For one of the $N=100$ simulated datasets, the pairs $(\beta^tX_i,{\hat{b}} X_i)$,
$i=1,\dots,n$ are represented on Figure~\ref{compasirtik} for Model 1 and on Figure~\ref{compasirtiknew} for Model 2. In order to make a comparison, the estimated axis ${\hat{b}}$ is computed with the classical SIR method and the PCA+Tikhonov method. It appears clearly that the SIR method leads to very bad oriented estimator of $\beta$ which is not the case with the PCA+Tikhonov method.

\begin{figure}[p]
\begin{center}
\epsfig{figure=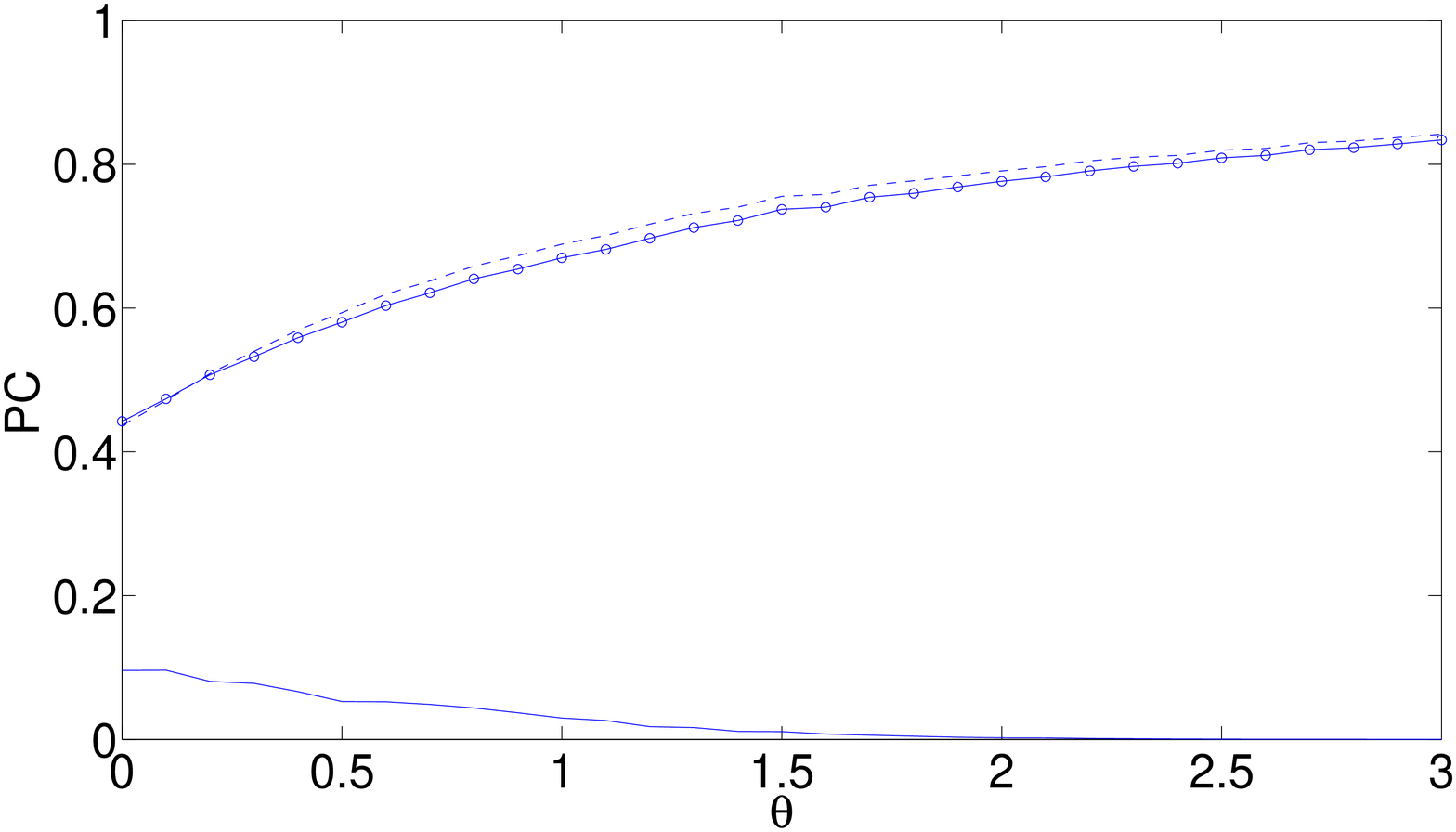,width=0.6\textwidth,height=0.5\textwidth}
\end{center}
\caption{Sensibility of GRSIR with respect to the condition number of
the covariance matrix for Model 1. Here, we use the optimal $\tau$.
Horizontally: $\theta$, vertically: the criterion MSC.
Continuous line: $\Omega_0$ (SIR), "-o-" line: $\Omega_1$ (ridge),
dashed line: $\Omega_3$ (Tikhonov).}
\label{figcondi}
\end{figure}

\begin{figure}[p]
\begin{center}
\epsfig{figure=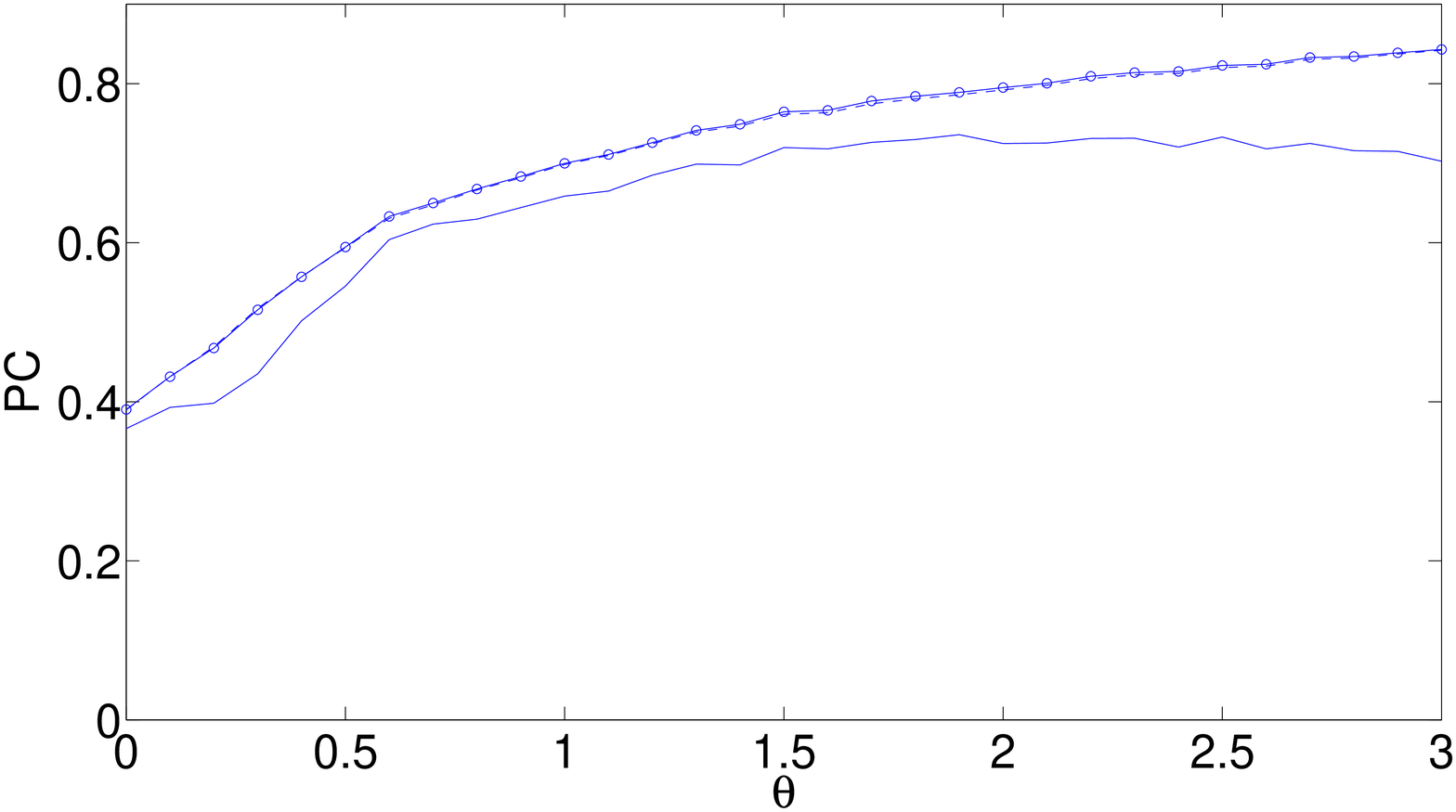,width=0.6\textwidth,height=0.5\textwidth}
\end{center}
\caption{Sensibility of GRSIR with respect to the condition number of
the covariance matrix for Model 1. The cut-off dimension is chosen to $d=20$ and we use the optimal $\tau$..
Horizontally: $\theta$, vertically: the criterion MSC.
Continuous line: $\Omega_2$ (PCA+SIR), "-o-" line: $\Omega_4$ (PCA+ridge),
dashed line: $\Omega_5$ (PCA+Tikhonov).}
\label{figacpcondi}
\end{figure}

Second experiment: The robustness with respect
to the condition number is investigated by varying $\theta$
in $\{0,0.1,0.2,\dots, 3\}$. For each value within this set,
the optimal regularization parameter has been selected for each method
and the corresponding MSC criterion is displayed on Figure~\ref{figcondi}
and Figure~\ref{figacpcondi}. This is done only for Model 1.
Clearly, the classical SIR method is very sensitive to the
ill-conditioning of the covariance matrix.
For all the other considered methods, results are getting better
while the condition number increases. 
Note that ridge and Tikhonov methods as well as
PCA+ridge and PCA+Tikhonov yield very similar results.

\begin{figure}[p]
\begin{center}
\epsfig{figure=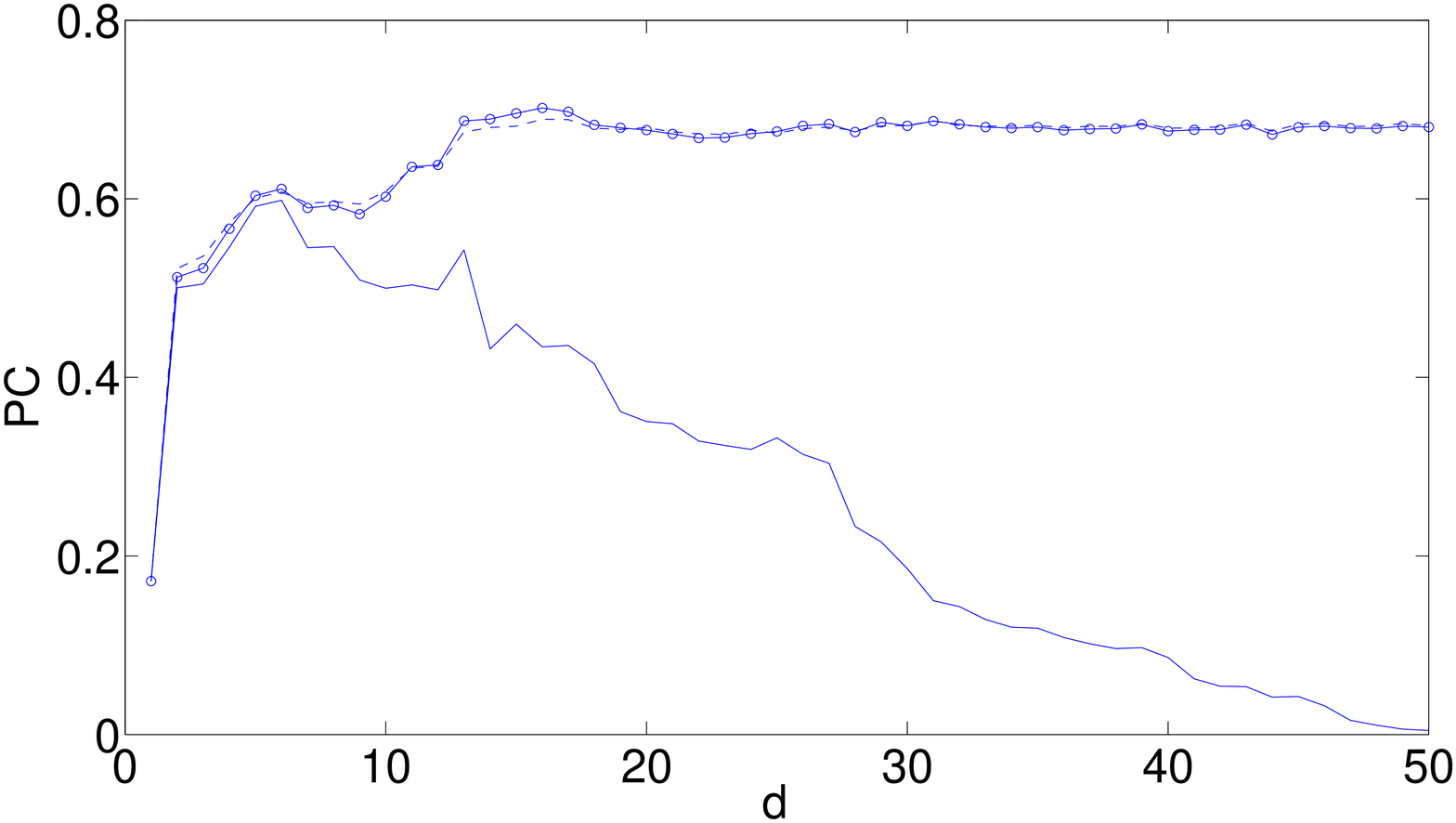,width=0.6\textwidth,height=0.5\textwidth}
\end{center}
\caption{Sensibility of GRSIR with respect to 
the cut-off dimension for Model 1.
Horizontally: $d$, vertically: the MSC criterion.
Continuous line: $\Omega_2$ (PCA+SIR), "-o-" line: $\Omega_4$ (PCA+ridge),
dashed line: $\Omega_5$ (PCA+Tikhonov). Here, $\theta=2$ and the optimal $\tau$ is used.}
\label{figbestd50}
\end{figure}

Third experiment: Illustration of the
role of $d$ in PCA+SIR, PCA+ridge and PCA+Tikhonov methods.
Here, the condition number is fixed by choosing $\theta=2$.
 For each value of $d$
in $\{0,1,\dots,p\}$, the optimal regularization parameter has been selected for each method
and the corresponding MSC criterion is displayed on Figure~\ref{figbestd50}. Only the Model 1 is considered here.
One can see that the PCA+SIR method is very sensitive to $d$.
Indeed, if $d$ is large, then this approach reduces to SIR,
whose accuracy is low for large dimensions.
At the opposite, PCA+ridge and PCA+Tikhonov results remain
stable as $d$ increases, since these methods get close to ridge
and Tikhonov methods respectively.

%%%%%%%%%%%%%%%%%%%%%%%%%%%%%%%%%%%%%%%%%%%%%%%%%%%%%%%%%%%%%%%%%%%%%%%
\section{Retrieval of Mars surface physical properties from hyperspectral images}
%%%%%%%%%%%%%%%%%%%%%%%%%%%%%%%%%%%%%%%%%%%%%%%%%%%%%%%%%%%%%%%%%%%%%%%

We propose here to apply GRSIR in the context of a nonlinear inverse problem in remote sensing. Hyperspectral remote sensing is a promising space technology regularly selected by agencies with regard to the exploration and observation of planets.
It allows to collect for each pixel of a scene, the intensity of light energy reflected from materials as it varies across different wavelengths. Hundreds spectels in the visible and near infra-red are recorded, making it possible to observe a continuous spectrum for each image cell. The analyze of these spectral signatures allows to identify the physical, chemical or mineralogical properties of the surface and of the atmosphere that may help to understand the geological and climatological history of planets.

Our goal is to evaluate the physical properties of surface materials on Mars planet from hyperspectral images collected by the OMEGA instrument aboard Mars express spacecraft. The used approach is based on the estimation of the functional relationship $G$ between some physical parameters $Y$ and observed spectra $X$. For this purpose, a database of synthetic spectra is generated by a physical radiative transfer model and used to estimate $G$. The high dimension of spectra ($p=184$ wavelengths) is reduced by using GRSIR. Results are compared with the traditional SIR approach.

\subsection{Data}

In this application, we focus on an observation of the south pole of Mars at the end of summer. It has been collected during orbit 61 by the French imaging spectrometer OMEGA on board Mars Express Mission. A detailed analysis of this image~\cite{doute} revealed that this portion of Mars mainly contains water ice, carbon dioxide and dust. This has led to the physical modeling of individual spectra with a surface reflectance model. This model allows the generation of 12.000 synthetic spectra with the corresponding parameters that constitute a learning database. Here, we focus on the terrain unit of strong CO$_2$ concentration determined by a classification method based on wavelets~\cite{schmidt}. It contains approximately 9,000 spectra to reverse. The 5 most important parameters characterizing the morphology of these spectra are the proportions of CO$_2$ ice, water ice and dust; and the grain sizes of water ice and CO$_2$ ice. In the sequel only one parameter, the grain size of CO$_2$ ice is presented. A detailed analysis for all other parameters can be found in~\cite{cbm}.

\subsection{Methodology}

\begin{figure}
\begin{center}
 \noindent\includegraphics[scale=0.30]{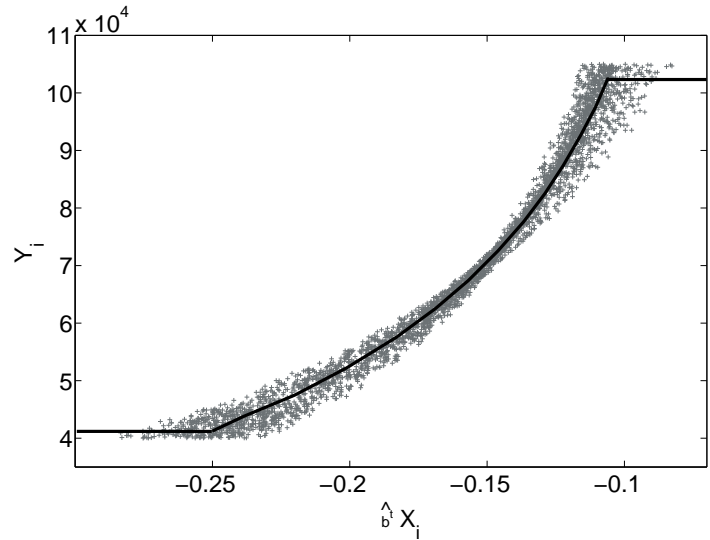}
\end{center}
 \caption{Functional relationship between reduced spectra $\hat \beta^{t} X$ on the first GRSIR direction and $Y$, the grain size of CO$_2$. Horizontally: reduced spectra from the learning database on the first GRSIR direction. Vertically: Grain size of CO$_2$.} \label{RelFn}
 \end{figure}

Two methods are compared in order to reverse the hyperspectral image: 
SIR and GRSIR with the PCA + ridge approach.
Both methods allow to choose a lower dimensional regressor space, sufficient enough to describe the relationship of interest. In practice, it appears  that a unidimensional regressor space gives satisfactory results~\cite{cbm}. Thus, we assume that there exists a function $g$ from $\R$ to $\R$ such that:
$G(X)=g(<\beta, X>)$.
As a consequence, the estimation of the relationship $G$ reduces to computing one direction $\beta$, the univariate function $g$ being estimated by piecewise linear regression.
In the PCA+ridge approach, the cut-off dimension is fixed to $d=20$.
The regularization parameter, fixed to $\tau=10^{-4.2}$, is chosen to minimize 
the mean squared error when estimating the grain size of CO$_2$ ice of the learning database itself.
An example of the relationship between reduced spectra and the grain size of CO$_2$ estimated using GRSIR is presented on Figure~\ref{RelFn}. It shows that the relationship is nonlinear and that one direction seems to be sufficient to estimate the grain size of CO$_2$.

\subsection{Results}

\begin{figure}
\begin{center}
 \noindent\includegraphics[scale=0.30]{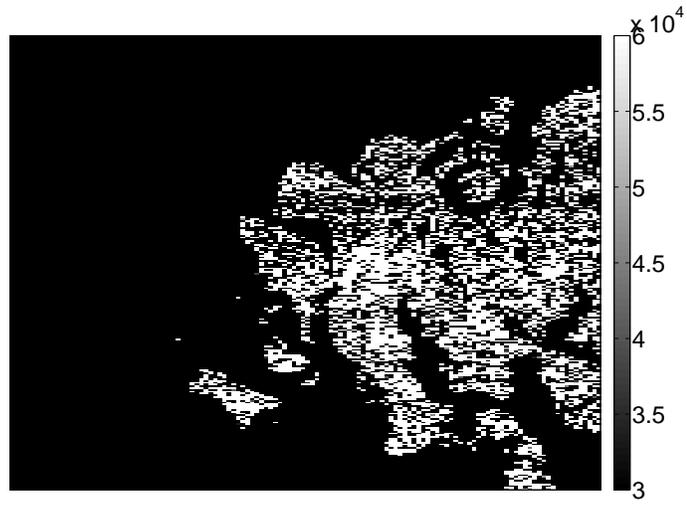}
 \caption{Grain size of CO$_2$ ice estimated by SIR on an hyperspectral image observed on Mars during orbit 61.} \label{SIR61}
\end{center}
 \end{figure}

\begin{figure}
\begin{center}
 \noindent\includegraphics[scale=0.30]{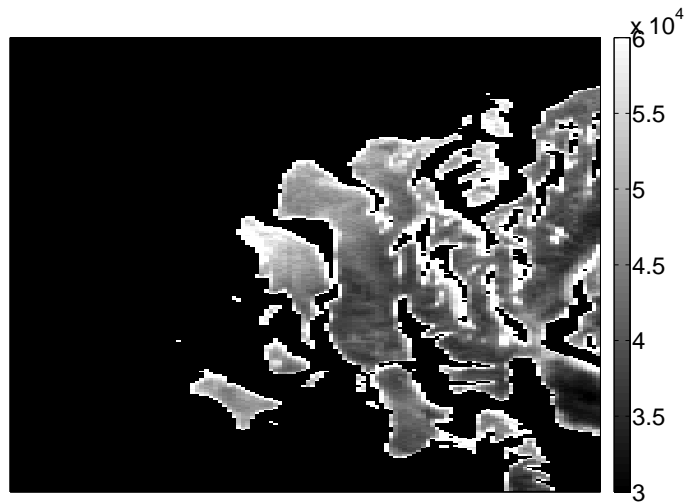}
 \caption{Grain size of CO$_2$ ice estimated by GRSIR (PCA+ridge) on an hyperspectral image observed on Mars during orbit 61.} \label{GRSIR61}
\end{center}
 \end{figure}

With GRSIR methodology, the inversion of the image from orbit 61 shows a smooth mapping of the grain size of CO$_2$ (see Figure~\ref{GRSIR61}) making it possible to distinguish some areas with great sizes of CO$_2$ ice on the boundaries and some areas with small values inside the cap. On the opposite, with the traditional SIR approach, estimates assume a small number of values that seem to be distributed randomly and correspond to the minimum and maximum values of the parameter in the learning database (see Figure~\ref{SIR61}). 
These poor results can be explained by the very high condition number, about
$10^{14}$, of the empirical covariance matrix.
Since no ground data is available, it is quite difficult to quantify the accuracy of GRSIR estimates. However, comparisons with other approaches or with estimates from other hyperspectral images from the same portion of Mars~\cite{cbm} give consistent results. GRSIR approach then appears promising for model inversion in remote sensing.

%%%%%%%%%%%%%%%%%%%%%%%%%%%%%%%%%%%%%%%%%%%%%%%%%%%%%%%%%%%%%%%%%%%%%%%
\section{Concluding remarks}
%%%%%%%%%%%%%%%%%%%%%%%%%%%%%%%%%%%%%%%%%%%%%%%%%%%%%%%%%%%%%%%%%%%%%%%

A new framework has been presented to regularize the SIR method.
It provides new interpretations of the existing methods as well as
the construction of new regularization techniques. 
Among them, it appears that the PCA+ridge and PCA+Tikhonov
methods are interesting alternatives to the existing PCA+SIR~\cite{Chiar,LiLi}
and ridge~\cite{Zhong} methods. The use of a "cut-off" dimension $d$
permits to limit the sensitivity to the choice of the regularization
parameter $\tau$. The choice of this dimension itself seems not crucial
since for large value the above methods are close to the ridge
and Tikhonov techniques.
In our experiments, the choice $d=p/2$ appears as a good heuristics
in most situations.
Of course, an automatic selection of $(\tau,d)$ would be of interest.
To this end, the construction of a generalized cross-validation criterion
is under investigation.
We also plan to study the introduction on non-Gaussian priors
in order to obtain $L_1-$ penalizations and thus sparse estimates of $\beta$.

%%%%%%%%%%%%%%%%%%%%%%%%%%%%%%%%%%%%%%%%%%%%%%%%%%%%%%%%%%%%%%%%%%%%%%%

%%%%%%%%%%%%%%%%%%%%%%%%%%%%%%%%%%%%%%%%%%%%%%%%%%%%%%%%%%%%%%%%%%%%%%%
\section*{Acknowledgments}
This research is partially supported by IAP research network grant nr. P6/03 of the Belgian government (Belgian Science Policy). 
We are grateful to Sylvain Dout\'e for his important contribution to the data
and to the referees whose remarks led to an important
improvement of the presentation.

%%%%%%%%%%%%%%%%%%%%%%%%%%%%%%%%%%%%%%%%%%%%%%%%%%%%%%%%%%%%%%%%%%%%%%%

%%%%%%%%%%%%%%%%%%%%%%%%%%%%%%%%%%%%%%%%%%%%%%%%%%%%%%%%%%%%%%%%%%%%%%%
\section*{Appendix: Proofs}
%%%%%%%%%%%%%%%%%%%%%%%%%%%%%%%%%%%%%%%%%%%%%%%%%%%%%%%%%%%%%%%%%%%%%%%

\paragraph{Proof of Lemma~\ref{lemmaG}.}

Let us remark that 
\begin{equation}
\label{eqG}
\Delta:=G(\mu,V,b,c)-\log\det V = \frac{1}{n} \sum_{i=1}^n Z_i^t V^{-1} Z_i
\end{equation}
where we have defined for $i=1,\dots,n$,
\begin{eqnarray*}
Z_i&=& \mu+ s^t(Y_i) c Vb-{X}_i \\
&=& (\mu-\bar{X} + \bar{s}^t  c Vb )+ ( s(Y_i) - \bar{s} )^t  c Vb - (X_i-\bar{X})\\
&:=& Z_{1,i} + Z_{2,i} - Z_{3,i}.
\end{eqnarray*}
Since $Z_{2}$ and $Z_{3}$ are centered, replacing the previous
expansion in~(\ref{eqG}) yields
$$
\Delta= \frac{1}{n}  \sum_{i=1}^n Z_{1,i}^t V^{-1} Z_{1,i} + 
 \frac{1}{n}\sum_{i=1}^n Z_{2,i}^t V^{-1} Z_{2,i}+
 \frac{1}{n}\sum_{i=1}^n Z_{3,i}^t V^{-1} Z_{3,i}-
 \frac{2}{n}\sum_{i=1}^n Z_{2,i}^t V^{-1} Z_{3,i},
$$
where
\begin{eqnarray*}
 \frac{1}{n}  \sum_{i=1}^n Z_{1,i}^t V^{-1} Z_{1,i}  &=&(\mu-\bar{X} + \bar{s}^t  c Vb )^t V^{-1} (\mu-\bar{X} + \bar{s}^t  c Vb )\\
\frac{1}{n}\sum_{i=1}^n Z_{2,i}^t V^{-1} Z_{2,i} &=& (c^t W c)(b^t V b)\\
\frac{1}{n}\sum_{i=1}^n Z_{3,i}^t V^{-1} Z_{3,i} &=& 
\frac{1}{n}\sum_{i=1}^n \mbox{tr}((X_i-\bar{X})^t V^{-1} (X_i-\bar{X}))
= \frac{1}{n}\sum_{i=1}^n \mbox{tr}( V^{-1} (X_i-\bar{X})(X_i-\bar{X})^t)\\
&=& \mbox{tr}( V^{-1} \hat\Sigma)\\
 \frac{2}{n}\sum_{i=1}^n Z_{2,i}^t V^{-1} Z_{3,i} &=&
2c^t Mb,
\end{eqnarray*}
and the conclusion follows.\CQFD

%%%%%%%%%%%%%%%%%%%%%%%%%%%%%%%%%%%%%%%%%%%%%%%%%%%%%%%%%%%%%%%%%%%%%% 
\paragraph{Proof of Proposition~\ref{propproj}.}

Annulling the gradients of ${G}(\mu,V,b,c)$ yields the system of equations
\begin{eqnarray}
\label{eq0tmp}
\frac{1}{2}\nabla_\mu G&=& {\hat V}^{-1}(\hat\mu-\bar{X}+ \bar{s}^t \hat c \hat  V\hat  b) = 0, \\
\label{eq1tmp}
\frac{1}{2}\nabla_b G&=& \left(({\hat c}^t \bar{s})^2+{\hat c}^t W \hat c\right) \hat  V\hat  b - M^t \hat c + 
({\hat c}^t \bar{s})(\hat \mu-\bar{X})= 0, \\
\label{eq2tmp}
\frac{1}{2}\nabla_c G&=& ({\hat b}^t \hat V\hat  b)(\bar{s}\bar{s}^t+  W) \hat c - M\hat  b +
(\hat \mu-\bar{X})^t \hat b \bar{s} = 0,\\ 
\label{eq3tmp}
\nabla_V G &=& {\hat V}^{-1} + \hat b {\hat b}^t \left(({\hat c}^t \bar{s})^2+{\hat c}^t W \hat c\right) 
-  {\hat V}^{-1} \left( (\hat\mu-\bar X)(\hat \mu-\bar X)^t +\hat \Sigma \right) {\hat V}^{-1}=0.
\end{eqnarray}
From (\ref{eq0tmp}), we have 
\begin{equation}
\label{eqhatmu}
\hat\mu=\bar{X}- \bar{s}^t \hat c \hat  V \hat b.
\end{equation}
Replacing
in~(\ref{eq1tmp}) and~(\ref{eq2tmp}) yields the simplified system of equations
\begin{eqnarray}
\label{eq1}
({\hat c}^t W \hat c) \hat  V \hat b &=& M^t \hat c, \\
\label{eq2}
({\hat b}^t \hat V \hat b)  W \hat c &=& M \hat b . 
\end{eqnarray}
Assuming $W$ is regular, equation~(\ref{eq2}) entails
$$
\hat c = W^{-1} M \hat b / {\hat b}^t \hat V \hat b
$$
and replacing in~(\ref{eq1}) yields
\begin{equation}
\label{eq3}
({\hat c}^t W \hat c) ({\hat b}^t \hat V \hat b)  \hat V \hat b =  M^t W^{-1} M \hat b.
\end{equation}
Now, multiplying (\ref{eq3tmp}) on the left and on the right by $\hat V$
and taking account of~(\ref{eqhatmu}) entail
\begin{equation}
\label{eqV}
\hat \Sigma = \hat V + ({\hat c}^t W \hat c) \hat V \hat b {\hat b}^t \hat V.
\end{equation}
As a consequence of (\ref{eqV}), we have
\begin{equation}
\label{relation}
\hat\Sigma \hat b = \left( 1 + ({\hat c}^t W \hat c) ({\hat b}^t \hat V \hat b) 
\right) \hat V \hat b,
\end{equation}
which means that $\hat\Sigma\hat b$ is proportional to $\hat V\hat b$.
Consequently, one also has
\begin{equation}
\label{eqprop}
\hat V \hat b = \theta(\hat b)\hat\Sigma \hat b,
\end{equation}
where we have defined
$$
\theta(\hat b)=\frac{{\hat b}^t \hat V  \hat b}{ {\hat b}^t \hat \Sigma  \hat b}.
$$
Substituting (\ref{eqprop}) in (\ref{eq3}) yields
$$
({\hat c}^t W \hat c) ({\hat b}^t \hat V \hat b)  \theta(\hat b) \hat \Sigma \hat b =  M^t W^{-1} M \hat b
$$
and thus $b$ is an eigenvector of $\hat\Sigma^{-1} M^t W^{-1} M$. Let us denote
by $\hat\lambda$ the associated eigenvalue 
\begin{equation}
\label{eqlambda}
\hat\lambda= ({\hat c}^t W \hat c) ({\hat b}^t \hat V \hat b)  \theta(\hat b). 
\end{equation}
Collecting (\ref{relation}) and (\ref{eqprop}) yields
$$
\frac{1}{\theta(\hat b)}=1 + ({\hat c}^t W \hat c) ({\hat b}^t \hat V \hat b),
$$
and thus the eigenvalue can be rewritten as
$$
\hat \lambda = 1 -  \theta(\hat b). 
$$
Moreover, we have, from (\ref{eq1}),
\begin{eqnarray}
\label{eqq1}
{\hat{c}}^t M \hat b &=&  \hat \lambda / \theta(\hat b), \\
\label{eqq2}
\mbox{tr}(\hat\Sigma{\hat V}^{-1})&=&p +  \hat \lambda / \theta(\hat b),  \\
\nonumber
\log\det\hat V &=& \log\det\hat\Sigma - \log\det\left( I_p+({\hat c}^t W \hat c) \hat V\hat b {\hat b}^t\right)\\
\label{eqq3}
&=& \log\det\hat\Sigma - \log\left(1+ \hat \lambda / \theta(\hat b) \right),
\end{eqnarray}
entailing
\begin{eqnarray*}
 G(\hat\mu,\hat V,\hat b,\hat c) &=& p+ \log\det\hat\Sigma - \log\left(1+  \hat \lambda / \theta(\hat b) \right) 
\\
&=&  p+ \log\det\hat\Sigma + \log(1-\hat\lambda).
\end{eqnarray*}
As a consequence, to minimize $G$, $\hat\lambda$ should be the largest
eigenvalue.  Finally, let us consider
$$
V_0= \hat\Sigma -\frac{\hat\lambda}{{\hat b}^t\hat\Sigma \hat b} \hat\Sigma \hat b {\hat b}^t \hat \Sigma,
$$
leading to 
\begin{eqnarray*}
V_0 + ({\hat c}^t W \hat c) V_0 \hat b {\hat b}^t V_0 &=&
\hat\Sigma + \left( ({\hat c}^t W \hat c) \theta^2(\hat b) - 
\frac{\hat\lambda}{{\hat b}^t\hat\Sigma \hat b}\right)
 \hat\Sigma \hat b {\hat b}^t \hat \Sigma
 \\
&=& \hat\Sigma + \frac{\hat\lambda  \theta({\hat b})}{{\hat b}^t V_0 \hat b}
 \left( ({\hat c}^t W \hat c) ({\hat b}^t V_0 \hat b)
 \theta(\hat b) - 
\hat\lambda \right)
 \hat\Sigma \hat b {\hat b}^t \hat \Sigma\\
&=& \hat\Sigma,
\end{eqnarray*}
in view of (\ref{eqlambda}) and thus $V_0$ verifies equation~(\ref{eqV}).  \CQFD

\paragraph{Proof of Corollary~\ref{corodiscr}.} 
%%%%%%%%%%%%%%%%%%%%%%%%%%%%%%%%%%%%%%%%%%%%%%%%%%%%%%%%%%%%%%%%%%%%%% 
Let us remark that, under~(\ref{baseindic}), the coefficients of $W_{ij}$ of $W$
are given by $W_{ij}=f_i \I\{i=j\}-f_i f_j$ for all $(i,j)\in\{1,\dots,h\}^2$.
The inverse matrix of $W$ is 
$$
W^{-1}=\mbox{diag}\left(\frac{1}{f_1},\dots,\frac{1}{f_h}\right)+ \frac{1}{f_{h+1}}U,
$$
where $U$ is the $h\times h$ matrix defined by $U_{ij}=1$
for all $(i,j)\in\{1,\dots,h\}\times\{1,\dots,h\}$.
Since the $j$th line of $M$ is
given by $f_j (\bar{X}_j-\bar X)^t$ for all $j=1,\dots,h$
and taking account of $U^2=f_{h+1} U$,
we have
\begin{eqnarray}
\nonumber
M^t W^{-1} M &=& \sum_{j=1}^h f_j  (\bar{X}_j-\bar X)(\bar{X}_j-\bar X)^t +
\frac{1}{f_{h+1}} M^t U M\\
\label{eq123}
&=&  \sum_{j=1}^h f_j  (\bar{X}_j-\bar X)(\bar{X}_j-\bar X)^t +
\frac{1}{h f_{h+1}}  (M^t U) (M^t U)^t.
\end{eqnarray}
Now, remarking that all the columns of $M^t U$ are equal to
\begin{eqnarray*}
 \sum_{j=1}^h f_j  (\bar{X}_j-\bar X) &=&  \sum_{j=1}^h f_j  (\bar{X}_j-\bar X)  - f_{h+1} (\bar{X}_{h+1}-\bar X)\\
&=& - f_{h+1} (\bar{X}_{h+1}-\bar X),
\end{eqnarray*}
it follows that
$$
 (M^t U) (M^t U)^t = h f_{h+1}^2  (\bar{X}_{h+1}-\bar X) (\bar{X}_{h+1}-\bar X)^t
$$
and thus replacing in~(\ref{eq123}) yields
\begin{eqnarray*}
M^t W^{-1} M &=&  \sum_{j=1}^{h+1} f_j  (\bar{X}_j-\bar X)(\bar{X}_j-\bar X)^t\\
&=&\hat\Gamma.
\end{eqnarray*}
The result is then a consequence of Proposition~\ref{propproj}.  \CQFD

\paragraph{Proof of Lemma~\ref{lemmap}.} 
%%%%%%%%%%%%%%%%%%%%%%%%%%%%%%%%%%%%%%%%%%%%%%%%%%%%%%%%%%%%%%%%%%%%%% 
The joint distribution of $(X,b)$ given $Y$ denoted by $p(X,b|Y)$ 
is calculated as the product $p(X|Y,b) p(b|Y)$ where $p(X|Y,b)$ is given
by~(\ref{IRMproj}) and $p(b|Y)$ is given by~(\ref{apriori}).
The estimators are obtained by minimizing
\begin{eqnarray*}
J_{\Omega}(\mu,V,b,c)&=&- \frac{2}{n}\sum_{i=1}^n \log p(X_i,b|Y_i)
  =- \frac{2}{n}\sum_{i=1}^n \log p(X_i|Y_i,b)
  - \frac{2}{n}\sum_{i=1}^n \log p(b|Y_i)\\
&=& G(\mu,V,b,c) +  \frac{b^t\Omega^{-1}b}{1+\rho} \frac{1}{n} \sum_{i=1}^n ((s(Y_i)-\bar{s})^t c)^2 + C^{ste}\\
&=& G(\mu,V,b,c) + \frac{b^t V b}{b^t\Sigma b} (b^t\Omega^{-1}b) (c^tW c) + C^{ste},
\end{eqnarray*}
which is $G_{\Omega}(\mu,V,b,c)$ up to the constant $C^{ste}$.\CQFD

\paragraph{Proof of Proposition~\ref{propregul}.}
%%%%%%%%%%%%%%%%%%%%%%%%%%%%%%%%%%%%%%%%%%%%%%%%%%%%%%%%%%%%%%%%%%%%%% 
The proof is similar to the one of Proposition~\ref{propproj}.
First, remark that equation (\ref{eq0tmp}) still holds and thus
$\hat\mu=\bar{X}- \bar{s}^t \hat c \hat  V \hat b$.
Let us recall the following definitions
$$
\theta(\hat b)=\frac{{\hat b}^t \hat V  \hat b}{ {\hat b}^t \hat \Sigma  \hat b}\;\mbox{ and }\;
\eta(\hat b)=\frac{{\hat b}^t \Omega^{-1} \hat b}{ {\hat b}^t \hat \Sigma  \hat b}.
$$
Annulling the gradients of $G_\Omega(\mu,V,b,c)$ yields the system of equations
\begin{eqnarray}
\label{eq1tmpbis}
({\hat c}^t W \hat c) \left( \Omega \hat V \hat b + \theta(\hat b) \hat b
+ \eta({\hat b})
\Omega ( \hat V\hat b - \theta(\hat b) \hat\Sigma \hat b)\right)
 &=& \Omega M^t \hat c, \\
\label{eq2tmpbis}
({\hat b}^t \hat V\hat  b)(1+ \eta({\hat b})) W \hat c &=& M\hat  b,\\ 
\label{eq3tmpbis}
{\hat V}^{-1} + \hat b {\hat b}^t ({\hat c}^t W \hat c)
(1+ \eta({\hat b})) &=& {\hat V}^{-1} \hat \Sigma  {\hat V}^{-1}.
\end{eqnarray}
Assuming $W$ is regular, equation~(\ref{eq2tmpbis}) entails
$$
\hat c = W^{-1} M \hat b / (( 1+\eta(\hat b))({\hat b}^t \hat V \hat b))
$$
and replacing in~(\ref{eq1tmpbis}) yields
\begin{equation}
\label{eq333}
({\hat c}^t W \hat c) ({\hat b}^t \hat V \hat b) ( 1+\eta(\hat b)) 
\left( \Omega \hat V \hat b + \theta(\hat b) \hat b
+ \eta({\hat b})
\Omega ( \hat V\hat b - \theta(\hat b) \hat\Sigma \hat b)\right)
 =  \Omega M^t W^{-1} M \hat b.
\end{equation}
Now, multiplying (\ref{eq3tmpbis}) on the left and on the right by $\hat V$,
it follows that
$$
\hat \Sigma = \hat V + ({\hat c}^t W \hat c)( 1+\eta(\hat b))  \hat V \hat b {\hat b}^t \hat V,
$$
leading to
\begin{equation}
\label{comp1}
\hat\Sigma \hat b =\left( 1 + ({\hat c}^t W \hat c) ({\hat b}^t \hat V \hat b)( 1+\eta(\hat b))  \right) \hat V \hat b, 
\end{equation}
which means that $\hat\Sigma \hat b$ is proportional to $\hat V \hat b$.
As a consequence, one also has
\begin{equation}
\label{comp2}
\hat V \hat b = \theta(\hat b)\hat\Sigma \hat b,
\end{equation}
and replacing in~(\ref{eq333}) yields
$$
({\hat c}^t W \hat c) ({\hat b}^t \hat V \hat b) ( 1+\eta(\hat b)) \theta(\hat b) \left( \Omega \hat \Sigma  + I_p\right) \hat b
= \Omega M^t W^{-1} M \hat b,
$$
and thus $b$ is an eigenvector of $(\Omega\hat\Sigma+ I_p)^{-1}\Omega M^t W^{-1} M$. 
Let us denote by $\hat\lambda$ the associated eigenvalue 
$$
\hat\lambda= ({\hat c}^t W \hat c) ({\hat b}^t \hat V \hat b) ( 1+\eta(\hat b)) \theta(\hat b).
$$
Collecting (\ref{comp1}) and (\ref{comp2}) yields
$$
\frac{1}{\theta(\hat b)} =  1 + ({\hat c}^t W \hat c) ({\hat b}^t \hat V \hat b)( 1+\eta(\hat b))  
$$
and consequently the eigenvalue can be rewritten as
$$
\hat\lambda = 1 - \theta(\hat b).
$$
Now, let us remark that~(\ref{eqq1}), (\ref{eqq2}) and (\ref{eqq3})
still hold in this context
entailing
$$
G_\Omega(\hat\mu,\hat V,\hat b,\hat c) 
=  p+ \log\det\hat\Sigma + \log(1-\hat\lambda).
$$
As a consequence, to minimize $G_\Omega$, $\hat\lambda$ should be the largest
eigenvalue. 
 The end of the proof follows the same lines as the
one of Proposition~\ref{propproj}.  \CQFD

\paragraph{Proof of Proposition~\ref{propacpsir}.}
%%%%%%%%%%%%%%%%%%%%%%%%%%%%%%%%%%%%%%%%%%%%%%%%%%%%%%%%%%%%%%%%%%%%%% 
Let $P$ be the projection matrix on $S_d$:
$$
P=\sum_{j=1}^d q_j q_j^t,
$$
and, for all $i=1,\dots,n$ consider the projected predictor defined by
$\tilde X_i=P X_i$.
Introducing $\tilde \Gamma= P\hat\Gamma P$ the empirical
"between slices" matrix associated to $\tilde X_1,\dots,\tilde X_n$
and $\tilde\Sigma=P\hat\Sigma P$ the corresponding covariance matrix,
the PCA+SIR method finds $\tilde b$ such that
$$
\tilde\Gamma \tilde b = \tilde\lambda\tilde\Sigma \tilde b,
$$
where $\tilde\lambda\in\R$, or equivalently,
$$
P\hat\Gamma P \tilde b = \tilde \lambda P\hat\Sigma P \tilde b.
$$
Remarking that $P\hat\Sigma P=\hat\Sigma P$, we have, for all $\tau>0$,
$$
P\hat\Gamma P \tilde b = \frac{\tilde \lambda}{1+\tau} (P\hat\Sigma +\tau\hat\Sigma) P \tilde b,
$$
and defining $b=P\tilde b$ and $\lambda=\tilde\lambda/(1+\tau)$, it follows that
$$
P\hat\Gamma  b = \lambda (P\hat\Sigma +\tau\hat\Sigma) b.
$$
Since $P=\tau \hat\Sigma\Omega_2$, we thus have
$$
\hat\Sigma\Omega_2\hat\Gamma  b = \lambda (\hat\Sigma\Omega_2\hat\Sigma +\hat\Sigma)b,
$$
which means that $b$ is an eigenvector of
$(\Omega_2\hat\Sigma +I_p)^{-1} \Omega_2\hat\Gamma$.
Corollary~\ref{cororegul} concludes the proof,
{\it i.e.} $b$ is the GRSIR estimator with prior covariance matrix $\Omega_2$.\CQFD

\paragraph{Proof of Proposition~\ref{propreducdim}.}
%%%%%%%%%%%%%%%%%%%%%%%%%%%%%%%%%%%%%%%%%%%%%%%%%%%%%%%%%%%%%%%%%%%%%% 
Here, we adopt the notations introduced in the proof of
Proposition~\ref{propacpsir}. The GRSIR estimator $\tilde{b}$ 
computed on the
projected predictors $\tilde X_1,\dots,\tilde X_n$ verifies
$$
\tilde\Omega(\varphi)\tilde\Gamma\tilde b = \tilde\lambda
(\tilde\Omega(\varphi)\tilde\Sigma + I_p) \tilde b,
$$
or equivalently
$$
\tilde\Omega(\varphi)P\hat\Gamma P \tilde b = \tilde\lambda
(\tilde\Omega(\varphi) P \hat\Sigma P + I_p) \tilde b.
$$
Multiplying this equation by $P$ on the left, we obtain
$$
P\tilde\Omega(\varphi)P\hat\Gamma P \tilde b = \tilde\lambda
(P\tilde\Omega(\varphi) P \hat\Sigma P + P) \tilde b.
$$
Since $P\tilde\Omega(\varphi)P=\Omega(\varphi)$, and introducing
$b=P\tilde b$, it follows that
$$
\Omega(\varphi)\hat\Gamma  b = \tilde\lambda
(\Omega(\varphi)  \hat\Sigma  + I_p ) b,
$$
which means that $b$ is an eigenvector of
$(\Omega(\varphi)\hat\Sigma +I_p)^{-1} \Omega(\varphi)\hat\Gamma$.
Corollary~\ref{cororegul} concludes the proof,
{\it i.e.} $b$ is the GRSIR estimator with prior covariance matrix $\Omega(\varphi)$.\CQFD

%%%%%%%%%%%%%%%%%%%%%%%%%%%%%%%%%%%%%%%%%%%%%%%%%%%%%%%%%%%%%%%%%%%%%%

%%%%%%%%%%%%%%%%%%%%%%%%%%%%%%%%%%%%%%%%%%%%%%%%%%%%%%%%%%%%%%%%%%%%%%

%%%%%%%%%%%%%%%%%%%%%%%%%%%%%%%%%%%%%%%%%%%%%%%%%%%%%%%%%%%%%%%%%%%%%%

%%%%%%%%%%%%%%%%%%%%%%%%%%%%%%%%%%%%%%%%%%%%%%%%%%%%%%%%%%%%%%%%%%%%%%%%%%%

%%%%%%%%%%%%%%%%%%%%%%%%%%%%%%%%%%%%%%%%%%%%%%%%%%%%%%%%%%%%%%%%%%%%%%

%%%%%%%%%%%%%%%%%%%%%%%%%%%%%%%%%%%%%%%%%%%%%%%%%%%%%%%%%%%%%%%%%%%%%%

\newpage

\end{document}